\documentclass[openany, a4paper, 11pt]{article}
\usepackage{mathrsfs,verbatim}
\usepackage{amsfonts}
\usepackage{booktabs}
\usepackage{threeparttable}

 \usepackage{mathrsfs,amsfonts,amsmath}
 \usepackage{color,cases}
 \def\benumerate{\begin{enumerate}}\def\eenumerate{\end{enumerate}}
 \def\bitemize{\begin{itemize}}\def\eitemize{\end{itemize}}

 \def\beqlb{\begin{eqnarray}}
 \def\eeqlb{\end{eqnarray}}
 \def\beqnn{\begin{eqnarray*}}
 \def\eeqnn{\end{eqnarray*}}

 \def\proof{\noindent{\it Proof.~~}}
 \def\qed{\hfill$\Box$\medskip}

 \def\<{\langle}\def\>{\rangle}

 \def\mcr{\mathscr}\def\mbb{\mathbb}
 \def\mbf{\mathbf}\def\mrm{\mathrm}

 \def\ar{\!\!&}

 \def\d{\mrm{d}}\def\e{\mrm{e}}

  \newtheorem{theorem}{Theorem}[section]
 \newtheorem{definition}[theorem]{Definition}
 \newtheorem{lemma}[theorem]{Lemma}
 \newtheorem{corollary}[theorem]{Corollary}
 \newtheorem{proposition}[theorem]{Proposition}
 \newtheorem{remark}[theorem]{Remark}
 \newtheorem{condition}[theorem]{Condition}
 \newtheorem{example}{Example}[section]

 \def\blemma{\begin{lemma}}\def\elemma{\end{lemma}}
 \def\bproposition{\begin{proposition}}\def\eproposition{\end{proposition}}
 \def\btheorem{\begin{theorem}}\def\etheorem{\end{theorem}}
 \def\bcorollary{\begin{corollary}}\def\ecorollary{\end{corollary}}
 \def\bremark{\begin{remark}}\def\eremark{\end{remark}}
 \def\bcondition{\begin{condition}}\def\econdition{\end{condition}}

 \def\benumerate{\begin{enumerate}}\def\eenumerate{\end{enumerate}}
 \def\bitemize{\begin{itemize}}\def\eitemize{\end{itemize}}

 \def\beqlb{\begin{eqnarray}}\def\eeqlb{\end{eqnarray}}
 \def\beqnn{\begin{eqnarray*}}\def\eeqnn{\end{eqnarray*}}

 \def\ar{\!\!\!&}

 \def\mcr{\mathscr}\def\mbb{\mathbb}\def\mbf{\mathbf}\def\mrm{\mathrm}

 \def\proof{\noindent{\it Proof.~~}}\def\qed{\hfill$\Box$\medskip}

\def\d{\mrm{d}}

\begin{document}

\centerline{\Large\textbf{Mixed state branching evolution for cell division models}}

\bigskip

\centerline{Shukai Chen,\,\footnote{Supported by National Key R\&D Program of China (No. 2022YFA1006003) and China Postdoctoral Science Foundation (No. 2022M720735).} Lina Ji\,\footnote{Supported by NSFC grant (No. 12301167 and No. 12271029), Guangdong Basic and Applied Basic Research Foundation (No. 2022A1515110986) and Guangdong Young Innovative Talents Project (No. 2022KQNCX105).} and Jie Xiong\,\footnote{Supported by National Key R\&D Program of China grant (No. 2022YFA1006102) and NSFC grant (No. 11831010).}}

\medskip
\centerline{{\it School of Mathematics and Statistics} \& {\it Key Laboratory of Analytical Mathematics}}

\centerline{{\it and Applications
(Ministry of Education)} \&} 

\centerline{{\it Fujian Provincial Key Laboratory of Statistics and Artificial Intelligence} \&}

 \centerline{{\it Fujian Key
Laboratory of Analytical Mathematics and Applications (FJKLAMA)} \&}
 
\centerline{\it Center for Applied Mathematics of Fujian
Province (FJNU), Fujian Normal University}

\centerline{\it Fuzhou 350007, People's Republic of China.}

\centerline{\it Faculty of Computational Mathematics and Cybernetics, Shenzhen MSU-BIT University}

\centerline{\it Shenzhen 518172, People's Republic of China.}

\centerline{\it Department of Mathematics and National Center for Applied Mathematics (Shenzhen)}

\centerline{\it Southern University of Science and Technology, Shenzhen 518055,  China.}

\smallskip

\centerline{{\it E-mail: skchen@mail.bnu.edu.cn}, {\it jiln@smbu.edu.cn} and {\it xiongj@sustech.edu.cn}}
\bigskip

{\narrower{\narrower

\centerline{\textbf{Abstract}}
We prove a scaling limit theorem for two-type Galton-Waston branching processes with interaction. The limit theorem gives rise to a class of mixed state branching processes with interaction using to simulate the evolution for cell division affected by parasites.  Such process can also be obtained by the pathwise unique solution to a stochastic equation system. Moreover, we present sufficient conditions for extinction with probability one and the exponential ergodicity in the $L^1$-Wasserstein  distance of such process in some cases.

\medskip


\par}\par}

%%%%%%%%%%%%%%%%%%%%%%%%%%%%%%%%%%%%%%%%%%%%%%%%%%%%%

\section{Introduction}
Let $\mbb{N}=\{0,1,2,...\}$. We consider a continuous time model in $D=[0,\infty)\times\mbb{N}$ for cells and parasites, where the behavior of cell division is infected by parasites. Informally, the quantity of parasites $(X(t))_{t\geq0}$ in a cell evolves as a continuous state branching process. The cells divide in continuous time at a rate $h(x,y)$ which may depend on the quantity of parasites $x$ and cells $y$. This framework is general enough to be applied for the modelling of other structured populations, for instance, grass-rabbit models in \cite{JX22}.

Many studies have been conducted on branching within branching processes to study such population dynamics in continuous time. In \cite{OW20}, the evolution of parasites is modelled by a birth-death process, while  the cells split according to a Yule process.  \cite{BT11} allows the quantity of parasites in a cell following a Feller diffusion. A continuous state branching process with jumps is considered to model the quantity of parasites in a cell in \cite{MS20}. In particular, \cite{OW20,BT11,MS20} describe cell populations in a tree structure, in this way, the population of cells at some time may be represented by a random point measure and associated martingale problems can be established by choosing test functions appropriately. Instead of \cite{OW20,BT11,MS20}, in this paper we ignore the tree structure and mainly focus on a parasite-cell model from a macro point of view. More precisely, we use a stochastic equation system to describe the sample path of such models,
\beqnn
\begin{cases}
X(t)=X(0)-b\int_0^t X(s)\,\d s+\int_0^t\sqrt{2cX(s)}\,\d B(s)+\int_0^t\int_0^{X(s-)}\int_0^\infty \xi\,\tilde{M}(\d s,\d u,\d \xi),\\
Y(t)=Y(0)+\int_0^t\int_0^{Y(s-)}\int_0^{h(X(s-),Y(s-))}\int_{\mbb{N}} (\xi-1)\,N(\d s,\d u,\d r,\d \xi),
\end{cases}
\eeqnn
where $b \in \mbb{R}$ and $c\geq0$ are constants, $(B(t))_{t \ge 0}$ is a standard Brownian motion, $h(\cdot,\cdot)\in C(\mbb{R}_+^2)^+,$ here  $C(\mbb{R}_+^2)^+$ is the collection of continuous positive functions defined on $\mbb{R}_+^2$. Let $(\xi\wedge \xi^2)\,m(\d \xi)$ be a finite measure on $(0, \infty)$ and $(p_\xi:\xi\in\mbb{N})$ be an offspring distribution satisfying $\sum_\xi \xi p_\xi<\infty$. Without losing generality, we assume $p_1 = 0.$ The above $M(\d s, \d u, \d \xi)$ is a Poisson random measure on $(0, \infty)^3$ with intensity $\d s\d u m(\d \xi)$, and $\tilde{M}(\d s,\d u,\d \xi)=M(\d s,\d u,\d \xi)-\d s\d um(\d \xi)$. The above $N$ is a Poisson random measure on $(0, \infty)^3\times\mbb{N}$ with intensity $\d s\d u\d r n(\d \xi)$, where $n(\d \xi) = p_\xi\sharp(\d\xi)$ and $\sharp(\cdot)=\sum_j\delta_j(\cdot)$ is the counting measure on $\mbb{N}$. Those three random elements ($B(t),$ $M$ and $N$) are independent of each other. Apparently, $(X(t))_{t\geq0}$ is indeed a {\it continuous-state branching process} (CB-process), see \cite{DL06,DL12}. In particular, when $h(\cdot,\cdot)\equiv r>0$, $(Y(t))_{t\geq0}$ is a standard continuous time Markov branching process with branching rate $r>0$ and offspring $(p_\xi,\xi\in\mbb{N})$, in this case, the system can be seen as a particular case of mixed state branching processes, which has been studied in \cite{CL21}.

For simplicity, we introduce another Poisson random measure and write it again by $N$ on $[0,\infty)^3\times\mbb{N}^{-1}$, $\mbb{N}^{-1}=\mbb{N}\cup\{-1\}$ with characteristic measure
$n(\d\xi)=p'_\xi\sharp(\d\xi)$,  $p'_\xi=p_{\xi+1}$. Then we can rewrite the system by
\begin{numcases}{}
X(t)=X(0)-b\int_0^t X(s)\,\d s+\int_0^t\sqrt{2cX(s)}\,\d B(s)+\int_0^t\int_0^{X(s-)}\int_0^\infty \xi\,\tilde{M}(\d s,\d u,\d \xi),\label{1.1a}\\
Y(t)=Y(0)+\int_0^t\int_0^{Y(s-)}\int_0^{h(X(s-),Y(s-))}\int_{\mbb{N}^{-1}} \xi\,N(\d s,\d u,\d r,\d \xi)\label{1.1b}.
\end{numcases}
In the rest of the paper, we use the stochastic equation system \eqref{1.1a}--\eqref{1.1b} to describe the parasite-cell model. In the literature on the theory of branching processes, the rescaling (in time or state) approach plays a valuable
role in establishing the connection among those branching processes,
see \cite{KW71}, \cite{L06,L11}, \cite{M09}, \cite{CL21} and \cite{Liu22} and the references therein. To the best of our knowledge, limited work has been done in branching processes with interactions. This leads to the
first purpose of this paper, and the establishment of strong uniqueness of solution to \eqref{1.1a}--\eqref{1.1b}. For a sequence of two-type Galton--Watson processes with interactions $\{(x_k(n),y_k(n))_{n\in\mbb{N}}\}_{k\geq1},$ we prove that $(x_k(\lfloor\gamma_kt\rfloor)/k,y_k(\lfloor\gamma_kt\rfloor))_{t\geq0}$ converges
in distribution to the solution to (\ref{1.1a})--(\ref{1.1b}) as $k \rightarrow \infty$ under suitable conditions. The pathwise uniqueness of solution to \eqref{1.1a}--\eqref{1.1b} is also given.

In addition, the second purpose of this paper is to study several long time behaviors of such process and we mainly obtain the extinction behavior and exponential ergodicity in the $L^1$-Wasserstein distance in some cases. The result of extinction behavior is inspired by \cite{LYZ19}. Furthermore, ergodicity is the foundation for a wide class of limit theorems and long-time behavior for Markov processes. Due to the nonlinearity of function $h$, the semigroup transition of $(X,Y)$ is not explicit. We obtain the ergodic property by a coupling approach, which has been proved to be effective in the study of ergodicity of nonlinear case, see \cite{C05,LW19} and the references therein.

We now introduce some notation. Let $\e_\lambda(z) =\e^{-\langle \lambda, z\rangle}$ for any $\lambda = (\lambda_1, \lambda_2) \in \mbb{R}_+^2$ and $z=(x, y) \in D,$ where $\langle \lambda, z\rangle = \lambda_1 x + \lambda_2 y.$ We use $C_b(D)$ to denote the set of all bounded functions $(x, y) \mapsto f(x, y)$ on $D$ with $x \mapsto f(x, \cdot)$ continuous. Let $C_b^2(D)$ be the subset of $C_b(D)$ with continuous bounded derivatives up to $2$nd order on $x$. Let $C_0^2(D)$ be the subset of $C_b^2(D)$ vanishing at infinity, and $C_c^2(D)$ be the subset of $C_0^2(D)$ with compact support. Define $C_b(\mbb{R}_+^2)$ to be the collection of all bounded continuous functions on $\mbb{R}_+^2$, which is a subset of $C_b(D).$ Let $C_b^{2,1}(\mbb{R}_+^2)$ be the subset of $C_b(\mbb{R}_+^2)$ with continuous bounded derivatives up to $2$nd order on $x$ and continuous bounded derivatives up to first order on $y$. Then we have $C_b^{2,1}(\mbb{R}_+^2) \subset C_b^2(D).$ Let ${\mbb{D}}([0, \infty), D)$ denote the space of c\`{a}dl\`{a}g paths from $[0, \infty)$ to $D$ furnished with the Skorokhod topology. In the integrals, we make the convention that, for $a \le b \in \mbb{R},$
\beqnn
\int_a^b = \int_{(a, b]}\quad \text{and}\quad \int_a^\infty = \int_{(a, \infty)}.
\eeqnn

This paper is structured as follows. The existence by a scaling limit of a sequence of two-type Galton-Watson processes with interaction and pathwise uniqueness of solution to (\ref{1.1a})--(\ref{1.1b}) are given in Section 2. In Section 3 the extinction behavior of the system is studied. In Section 4, an exponential ergodic property is proved under some conditions.

\section{Existence and pathwise uniqueness of solution}

The generator $A$ of $(X(t),Y(t))_{t\geq0}$ satisfying \eqref{1.1a}--\eqref{1.1b} is determined by
\beqlb\label{generator}
Af(z)&=& x\Big[-bf'_x+cf''_{xx}+\int_0^\infty\{f(x+\xi,y)-f(x,y)-\xi f'_x\} m(\d \xi)\Big]\nonumber\\
&& + \gamma(x,y)\int_{\mbb{N}^{-1}}\Big\{f(x,y+\xi)-f(x,y)\Big\}n(\d \xi)
\eeqlb
for any $f\in C_b^{2}(D)$, where $z = (x, y) \in D$ and $\gamma(x, y) = h(x, y)y$. %By Theorems \ref{t1.1} and \ref{comparison}, the pathwise uniqueness and comparison property hold for the solution to \eqref{1.1a}-\eqref{1.1b} under Condition \ref{condition1}.
Then
\beqlb\label{Lez}
A\mrm{e}_\lambda(z)=\mrm{e}_\lambda(z)
\Big[x\phi_1(\lambda_1)+\gamma(x,y)\phi_2(\lambda_2)\Big],
\eeqlb
where
\beqlb
\phi_1(\lambda_1)&=& b\lambda_1+c\lambda_1^2+\int_0^\infty(\mrm{e}^{-\lambda_1\xi}-1+\lambda_1\xi)\,m(\d\xi),\label{phi}\\
\quad\phi_2(\lambda_2)&=& \int_{\mbb{N}^{-1}}(\mrm{e}^{-\lambda_2\xi}-1)\,n(\d\xi).\label{phi2}
\eeqlb

We first consider the case of $h \in C_b(\mbb{R}_+^2)^+.$ Given the initial value $(x(0),y(0))\in\mbb{N}\times\mbb{N}$, let $(x(n),y(n))_{n\ge 0}$ be a two-dimensional process defined by
\beqlb\label{discrete process}
x(n)=\sum_{j=1}^{x(n-1)}\alpha_{n-1,j},\quad\quad  y(n)=\sum_{j=1}^{y(n-1)}\beta_{n-1,j,\theta(x(n-1),y(n-1))},\quad  n\geq1,
\eeqlb
where $\{\alpha_{n,j}:n\geq 0,j\geq1\}$ are integer-valued i.i.d. random variables with offspring distribution $(w(i): i\in\mbb{N})$. %$\{\beta_{n,j,\theta(x,y)}:n\geq1,j\geq1,x,y\in\mbb{N}\}$ are random variables depending on function $\theta$.
Given $x,y \in \mbb{N}$, the above  $\{\beta_{n,j,\theta(x,y)}:n\geq 0,j\geq1\}$ are i.i.d. integer-valued random variables with offspring distribution $(v^{\theta(x,y)}(i):i\in\mbb{N})$ depending on the function $\theta$. Let $g_1$ and $g^{\theta(x,y)}_2$ be the generating functions of $(w(i):i\in\mbb{N})$ and $(v^{\theta(x,y)}(i):i\in\mbb{N})$, respectively. It is known that $(x(n), y(n))_ {n\geq0}$ is a Markov process and we call it {\it two-type Galton-Watson process with interaction}. Suppose that there exists a sequence of two-type Galton-Watson processes with interaction $\{(x_k(n),y_k(n))_{n\geq0}\}_{k\geq1}$ with parameters $(g_{k,1},g_{k,2}^{\theta_k(x,y)})$. Let $\{\gamma_k\}_{k\geq1}$ be a sequence of positive numbers with $\gamma_k \rightarrow \infty$ as $k \rightarrow \infty$.  For $(x, y) \in \mbb{N}\times\mbb{N},$ we introduce several functions on $\mbb{R}_+$ as below:
\beqnn
&&\bar{\Phi}_{k,1}(\lambda_1)=k\gamma_k\log\Big[1-(k\gamma_k)^{-1}\Phi_{k,1}(\lambda_1)\mrm{e}^{\lambda_1/k}\Big],\\
&&\Phi_{k,1}(\lambda_1)=k\gamma_k\Big[\mrm{e}^{-\lambda_1/k}-g_{k,1}(\mrm{e}^{-\lambda_1/k})\Big],\\
&&\bar{\Phi}_{k,2}^{\theta_k(x,y)}(\lambda_2)=\gamma_k\log\Big[1-\gamma_k^{-1}\Phi_{k,2}^{\theta_k(x,y)}(\lambda_2)\mrm{e}^{\lambda_2}\Big],\\
&&\Phi_{k,2}^{\theta_k(x,y)}(\lambda_2)=\gamma_k\Big[\mrm{e}^{-\lambda_2}-g_{k,2}^{\theta_k(x,y)}(\mrm{e}^{-\lambda_2})\Big].
\eeqnn
Let $E_k = \{0, k^{-1}, 2k^{-1}, \cdots\}$ for each $k\geq1$. For any $x \in \mbb{R}_+,$ we take $x_k := \lfloor kx\rfloor/k.$ Then $x_k \in E_k$ and $|x_k - x| \le 1/k$. Let $D_k:= E_k \times \mbb{N}$. Then $D_k$ is a subset of $D$. We define a continuous-time stochastic process taking values on $D_k$ as $(X_k(t), Y_k(t))_{t \ge 0}:=(x_k(\lfloor \gamma_k t\rfloor)/k,y_k(\lfloor \gamma_k t\rfloor))_{t\geq0}.$
Denote $Z_k(t) = (X_k(t), Y_k(t))$ to simplify the notation. In order to state our main results in this section, we first present the assumption taken throughout this section.

\bcondition\label{condition}
\begin{itemize}
	\item[\,]
	\item[(2.1.1)] The sequence $\{\Phi_{k,1}(\lambda_1)\}_{k \ge 1}$ is uniformly Lipschitz in $\lambda_1$ on each bounded interval, and converges to a continuous function as $k \rightarrow \infty$;
	\item[(2.1.2)] $\gamma_k [1- v_k^{\theta_k(kx_k,y)}(1)] \rightarrow h(x, y)$ uniformly in $(x, y) \in \mbb{R}_+\times\mbb{N}$ as $k \rightarrow \infty;$
	\item[(2.1.3)] $\frac{v_k^{\theta_k(kx_k, y)}(\xi)}{1- v_k^{\theta_k(kx_k,y)}(1)} \rightarrow p_\xi$ for $\xi \in \mbb{N}\backslash \{1\}$ uniformly in $(x, y) \in \mbb{R}_+\times\mbb{N}$ as $k \rightarrow \infty.$
	\item[(2.1.4)] The sequence $\{\Phi_{k,2}^{\theta_k(kx_k,y)}(\lambda_2)\}_{k \ge 1}$ is uniformly Lipschitz in $\lambda_2$ on each bounded interval, where the Lipshcitz coefficient is independent from $x, y$.
\end{itemize}
\econdition

By \cite[Proposition 2.5]{L20}, under Condition (2.1.1), $\Phi_{k,1}(\lambda_1)$ converges to a function with representation \eqref{phi} as $k \rightarrow \infty$, see also \cite{L06,L11}. Moreover, there exists a constant $K > 0$ such that
\beqlb\label{Phi1'}
\sup_k\Phi_{k,1}'(0+) = \sup_k \gamma_k[g_{k,1}'(1-) - 1] \le K.
\eeqlb

\begin{example}
	Let $\{p_\xi: \xi = 0, 1, \cdots\}$ be an offspring distribution with $p_1 = 0.$ Let
	\beqnn
	v_k^{\theta_k(kx_k, y)}(\xi) = p_\xi \gamma_k^{-1/2}\left(1 - e^{-\gamma_k^{-1/2}h(x_k, y)}\right)
	\eeqnn
	for any $\xi \in \mbb{N}\backslash \{1\}$ and
	\beqnn
	v_k^{\theta_k(kx_k, y)}(1) = 1 - \gamma_k^{-1/2}\left(1 - e^{-\gamma_k^{-1/2}h(x_k, y)}\right).
	\eeqnn
	Then we have
	\beqnn
	\Phi_{k,2}^{\theta_k(kx_k,y)}(\lambda_2) = \gamma_k^{1/2}\left(1 - e^{-\gamma_k^{-1/2}h(x_k, y)}\right)\left(\mrm{e}^{-\lambda_2}-g(\mrm{e}^{-\lambda_2})\right)
	\eeqnn
	with $g(\mrm{e}^{-\lambda_2}) = \sum_{\xi= 0}^\infty p_\xi e^{-\lambda_2\xi}.$ It is easy to check that the above satisfies Condition (2.1.2)--(2.1.4).
\end{example}

\subsection{Main results}
We have the following statements about a scaling limit theorem of mixed state branching processes with interactions, existence and pathwise uniqueness of solution to
\eqref{1.1a}--\eqref{1.1b}.

\btheorem\label{existence}
Suppose that Condition \ref{condition} holds and $h \in C_b(\mbb{R}_+^2)^+$. Let $Z_k(0)$ converge in distribution to $Z(0)$ as $k \rightarrow \infty$ with  $\sup_k \mbb{E}[X_k(0) + Y_k(0)] < \infty$. Then $\{(Z_k(t))_{ t \ge 0}\}_{k \ge 1}$ converges in distribution on ${\mbb{D}}(0, \infty), D)$ to $(Z(t))_{t \ge 0}$, which is a solution to \eqref{1.1a}--\eqref{1.1b}.
\etheorem

\btheorem\label{t1.1}
Assume that $h \in C_b(\mbb{R}_+^2)^+$. For any given initial value $(X(0),Y(0))\in D$, the pathwise uniqueness holds for (\ref{1.1a})--(\ref{1.1b}) on $D$.
\etheorem

\begin{remark}
For the classical branching process, the transition semigroup can be determined uniquely by the log-Laplace transform since branching property. Then the linear hull of $\{\e^{-\lambda x}: \lambda \ge 0\}$ is the core for generator of the process. Therefore, the scaling limit theorem can be obtained in the sense of convergence in distribution on the skorokhod space by \cite[p.226 and pp.233-234]{EK86}. We refer to \cite{KW71,L06} and references therein for details. However, the branching property  is invalid for our model since the interaction. It is not enough to give estimation only when $f = e_\lambda.$ Therefore, in the following proof of Theorem \ref{existence}, we give the existence of solution to the martingale problem by tightness for general $f,$ which implies the existence of solution to \eqref{1.1a}-\eqref{1.1b}.

As we know, there is a unique strong solution to \eqref{1.1a}, which is a CB-process. Based on this, in the proof of Theorem \ref{t1.1}, we then construct the pathwise unique solution to \eqref{1.1b} by path stitching method.

\end{remark}

One sees that there is a unique positive strong solution to \eqref{1.1a}-\eqref{1.1b} by Theorems \ref{existence} and \ref{t1.1} when
$h \in C_b(\mbb{R}_+^2)^+.$ However, the boundedness assumption of $h$ is removed in the following result.

\btheorem\label{texistence}
Suppose that $h \in C(\mbb{R}_+^2)^+$. Then there exists a unique positive strong solution to \eqref{1.1a}--\eqref{1.1b}.
\etheorem

\subsection{Proof of main results}

\bproposition\label{function converge}
Under Conditions (2.1.2)--(2.1.3),  $\mrm{e}^{\lambda_2}\Phi^{\theta_k(kx_k,y)}_{k,2}(\lambda_2)$ converges to $-h(x,y)\phi_2(\lambda_2)$ uniformly for $(x,y,\lambda_2)\in \mbb{R}_+\times\mbb{N}\times\mbb{R}_+$ as $k\rightarrow\infty$, where $h \in C_b(\mbb{R}_+^2)^+$ and $\phi_2(\lambda_2)$ is given by \eqref{phi2}.
\eproposition
\proof One can see that
\beqnn
\ar\ar\mrm{e}^{\lambda_2}\Phi^{\theta_k(kx_k,y)}_{k,2}(\lambda_2)\\
\ar\ar\quad=\gamma_k\Big[1-\mrm{e}^{\lambda_2}g^{\theta_k(kx_k,y)}_{k,2}(\mrm{e}^{-\lambda_2})\Big]\\
\ar\ar\quad=\gamma_k\Big[1-\mrm{e}^{\lambda_2}\sum_{j=0}^{\infty}\mrm{e}^{-\lambda_2 j}v_k^{\theta_k(kx_k,y)}(j)\Big]\cr
\ar\ar\quad=\gamma_k \sum_{j=0}^\infty (1 - e^{-\lambda_2(j-1)})v_k^{\theta_k(kx_k,y)}(j)\cr
\ar\ar\quad=\gamma_k\left[1- v_k^{\theta_k(kx_k,y)}(1)\right]\int_{\mbb{N}^{-1}\backslash\{0\}} (1-e^{-\lambda_2 \xi}) \rho_k^{\theta_k(kx_k,y)}(\d \xi),
\eeqnn
where
\beqnn
\rho^{\theta_k(kx_k,y)}_{k}(\d \xi) &=& \frac{1}{1- v_k^{\theta_k(kx_k,y)}(1)}\sum_{j=0}^\infty v_k^{\theta_k(kx_k,y)}(j)\delta_{j-1}(\d \xi)\cr
&=& \frac{v_k^{\theta_k(kx_k,y)}(\xi+1)}{1- v_k^{\theta_k(kx_k,y)}(1)}\sharp(\d \xi)
\eeqnn
for $\xi \in \mbb{N}_{-1}\backslash\{0\}$ with $\sharp(\d \xi)$ being the counting measure on $\mbb{N}_{-1}.$ The result follows from Conditions (2.1.2)--(2.1.3).\qed

For $\lambda = (\lambda_1, \lambda_2) \in \mbb{R}_+^2,$ we then have
\beqlb\label{eZk}
\e_\lambda(Z_k(t))
&=& \e_\lambda(Z_k(0)) +\sum_{i=1}^{\lfloor \gamma_k t\rfloor} \left[\e_\lambda\left(Z_k\left(\frac{i}{\gamma_k}\right)\right) - e_\lambda\left(Z_k\left(\frac{i-1}{\gamma_k}\right)\right)\right]\cr
&=& \e_\lambda(Z_k(0)) +\sum_{i=1}^{\lfloor \gamma_k t\rfloor} \gamma_k^{-1} A_k \e_\lambda\left(Z_k\left(\frac{i-1}{\gamma_k}\right)\right) + M_{k,\lambda}(t)\cr
&=&  \e_\lambda(Z_k(0)) + \int_0^{\lfloor\gamma_kt\rfloor/\gamma_k} A_k e_\lambda(Z_k(s))\d s + M_{k,\lambda}(t),
\eeqlb
where
\beqlb\label{M}
M_{k,\lambda}(t) &=& \sum_{i=1}^{\lfloor \gamma_k t\rfloor}\Bigg\{ \left[\e_\lambda\left(Z_k\left(\frac{i}{\gamma_k}\right)\right) - \e_\lambda\left(Z_k\left(\frac{i-1}{\gamma_k}\right)\right)\right]\cr
&&  - \mbb{E}\left[\e_\lambda\left(Z_k\left(\frac{i}{\gamma_k}\right)\right) - \e_\lambda\left(Z_k\left(\frac{i-1}{\gamma_k}\right)\right)\Bigg| \mcr{F}_{\frac{i-1}{\gamma_k}}\right]\Bigg\}
\eeqlb
is a martingale and for $z = (x, y) \in D,$
\beqnn
A_k\mrm{e}_\lambda(z)=\gamma_k\bigg[(g_{k,1}(\mrm{e}^{-\lambda_1/k}))^{ kx_k}\cdot
(g_{k,2}^{ \theta_k(kx_k,y)}(\e^{-\lambda_2}))^y-\mrm{e}_\lambda(z)\bigg].
\eeqnn
One can check that
\beqnn
A_k\mrm{e}_\lambda(z)
=\mrm{e}_\lambda(z)\Big[x\bar{\Phi}_{k,1}(\lambda_1)
+y\bar{\Phi}_{k,2}^{\theta_k(kx_k,y)}(\lambda_2)\Big]+o(1).
\eeqnn
By the above, \cite[Proposition 2.5]{L20} and Proposition \ref{function converge}, we have the following estimation.

\btheorem\label{elimit}
Suppose that Condition \ref{condition} holds. Then for any $\lambda > 0,$ we have
\beqnn
\lim_{k \rightarrow \infty}\sup_{z \in D_k}|A_k\e_\lambda(z) - A\e_\lambda(z)| = 0,
\eeqnn
where $A$ is the generator defined by \eqref{generator}.
\etheorem

\bproposition\label{moment}
Suppose that Condition \ref{condition} holds. Let $T > 0$ be a fixed constant and
$$
\sup\limits_k \mbb{E}[X_k(0) + Y_k(0)] < \infty.
$$ Then we have
\beqnn
\sup_k\sup_{0 \le t \le T}\mbb{E}[X_k(t) + Y_k(t)] < \infty.
\eeqnn
\eproposition
\proof
By \eqref{Phi1'} one sees that $0 \le g_{k,1}'(1-) \le K/\gamma_k + 1.$ Then for $t \in [\frac{i}{\gamma_k}, \frac{i+1}{\gamma_k}),$ we have
\beqnn
\mbb{E}[X_k(t)] &=&  k^{-1}\mbb{E}[x_k(\lfloor \gamma_k t\rfloor)]\cr
&=&   g_{k,1}'(1-)k^{-1}\mbb{E}[x_k(\lfloor \gamma_k t\rfloor - 1)]\cr
&\le &  ( K/\gamma_k + 1) k^{-1}\mbb{E}[x_k(\lfloor \gamma_k t\rfloor - 1)].
\eeqnn
By induction, we have $\mbb{E}[X_k(t)] \le ( K/\gamma_k + 1)^{\lfloor \gamma_kt\rfloor} \mbb{E}[X_k(0)].$ Moreover, by Condition (2.1.4), we have
\beqnn
\sup_k\left|\frac{\partial}{\partial \lambda_2}\Phi_{k,2}^{\theta_k(kx_k,y)}(\lambda_2)\Big|_{\lambda_2 = 0}\right| = \sup_k\gamma_k\left|\frac{\partial}{\partial z}g_{k,2}^{\theta_k(kx_k, y)}(z)\Big|_{z = 1} - 1\right|\le K.
\eeqnn
Similarly, for $t \in [\frac{i}{\gamma_k}, \frac{i+1}{\gamma_k}),$ one sees that
\beqnn
\mbb{E} [Y_k(t)] &=& \mbb{E}[y_k(\lfloor \gamma_kt\rfloor)] = \mbb{E}\left[\sum_{k = 1}^{y_k(n-1)}\mbb{E}\left[ \beta_{n-1,k,\theta_k(x_k(n - 1), y_k(n - 1))}\Big|x_k(n - 1), y_k(n - 1)\right]\right]\Bigg|_{n= \lfloor \gamma_kt\rfloor}\cr
&=& \mbb{E}\left[y_k(n-1)\cdot\frac{\partial}{\partial z}g_{k,2}^{\theta_k(x_k(n - 1), y_k(n - 1))}(z)\Bigg|_{z = 1}\right]\Bigg|_{n= \lfloor \gamma_kt\rfloor}\cr
&\le& (1 + K/\gamma_k) \mbb{E}[y_k(\lfloor \gamma_kt\rfloor - 1)].
\eeqnn
Then we get $\mbb{E}[Y_k(t)] \le ( K/\gamma_k + 1)^{\lfloor \gamma_kt\rfloor} \mbb{E}[Y_k(0)]$ by induction. The result follows.
\qed

Let $\{\tau_k: k \ge 1\}$ be a sequence of bounded stopping times, and $\{\delta_k: k \ge 1\}$ be a sequence of positive constants with $\delta_k \rightarrow 0$ as $k \rightarrow \infty.$ For a fixed constant $T > 0,$ we assume that
\beqnn
0 \le \tau_k \le \tau_k + \delta_k \le T.
\eeqnn
\bproposition\label{p1.3}
Suppose that Condition \ref{condition} holds and $h \in C_b(\mbb{R}_+^2)^+$. Then for any $\lambda \in \mbb{R}_+^2,$ we have
\beqnn
\lim_{k \rightarrow \infty} \mbb{E}\left[\left|\e_\lambda(Z_k(\tau_k + \delta_k)) - \e_\lambda(Z_k(\tau_k))\right|^2\right] = 0.
\eeqnn
\eproposition

\proof
For any $\lambda \in \mbb{R}_+^2,$ by \eqref{eZk} we have
\beqnn
\ar\ar\mbb{E}\left[\left|\e_\lambda(Z_k(\tau_k + \delta_k)) - \e_\lambda(Z_k(\tau_k))\right|^2\right]\cr \ar\ar\qquad\le \left|\mbb{E}\left[\e_{2\lambda}(Z_k(\tau_k + \delta_k)) - \e_{2\lambda}(Z_k(\tau_k))\right]\right|\cr
\ar\ar\qquad\quad + \left|\mbb{E}\left[2\e_\lambda(Z_k(\tau_k))[\e_\lambda(Z_k(\tau_k + \delta_k)) - \e_\lambda(Z_k(\tau_k))]\right]\right|\cr
\ar\ar\qquad\le I_1 + I_2 + I_3,
\eeqnn
where
\beqnn
I_1 &=&  \left|\mbb{E}\left[\int_{\lfloor\gamma_k\tau_k\rfloor/\gamma_k}^{\lfloor\gamma_k(\tau_k+\delta_k)\rfloor/\gamma_k} A_k\e_{2\lambda}(Z_k(s))\d s\right]\right|,\cr
I_2 &=& \left|\mbb{E}\left[2\e_\lambda(Z_k(\tau_k)) \int_{\lfloor\gamma_k\tau_k\rfloor/\gamma_k}^{\lfloor\gamma_k(\tau_k+\delta_k)\rfloor/\gamma_k} A_k\e_{\lambda}(Z_k(s))\d s\right]\right|,\cr
I_3 &=& \left|\mbb{E}\left[2\e_\lambda(Z_k(\tau_k)) \left(M_{k,\lambda}(\tau_k+\delta_k) - M_{k,\lambda}(\tau_k)\right) \right]\right|.
\eeqnn
Then by \eqref{Lez}, Theorem \ref{elimit} and Proposition \ref{moment}, one can see that
\beqnn
I_1
&\le&  \mbb{E}\left[\int_{\lfloor\gamma_k\tau_k\rfloor/\gamma_k}^{\lfloor\gamma_k(\tau_k+\delta_k)\rfloor/\gamma_k} \left|A_k \e_{2\lambda}(Z_k(s)) - A \e_{2\lambda}(Z_k(s)) \right|\d s\right]\cr
&& + \mbb{E} \left[\int_{\lfloor\gamma_k\tau_k\rfloor/\gamma_k}^{\lfloor\gamma_k(\tau_k+\delta_k)\rfloor/\gamma_k} \left|A \e_{2\lambda}(Z_k(s)) \right|\d s\right] \le K \delta_k.
\eeqnn
Similarly,
\beqnn
I_2 &\le& 2\mbb{E}\left[\int_{\lfloor\gamma_k\tau_k\rfloor/\gamma_k}^{\lfloor\gamma_k(\tau_k+\delta_k)\rfloor/\gamma_k} \left|A_k \e_{\lambda}(Z_k(s))\right| \d s\right]\cr
&\le& 2\mbb{E}\left[\int_{\lfloor\gamma_k\tau_k\rfloor/\gamma_k}^{\lfloor\gamma_k(\tau_k+\delta_k)\rfloor/\gamma_k} \left|A_k \e_{\lambda}(Z_k(s)) - A\e_{\lambda}(Z_k(s))\right| \d s\right]\cr
&& + 2\mbb{E}\left[\int_{\lfloor\gamma_k\tau_k\rfloor/\gamma_k}^{\lfloor\gamma_k(\tau_k+\delta_k)\rfloor/\gamma_k} \left|A \e_{\lambda}(Z_k(s))\right| \d s\right] \le K\delta_k.
\eeqnn
Moreover, recall that $h$ is bounded, it follows from \eqref{M} and Doob's stopping theorem that
\beqnn
&&\mbb{E}\left[\e_\lambda(Z_k(\tau_k)) \left(M_{k,\lambda}(\tau_k+\delta_k) - M_{k,\lambda}(\tau_k)\right) \right]\cr
&&\qquad= \mbb{E}\left[\mbb{E}\left[\e_\lambda(Z_k(\tau_k)) \left(M_{k,\lambda}(\tau_k+\delta_k) - M_{k,\lambda}(\tau_k)\right)|\mcr{F}_{\lfloor\gamma_k\tau_k\rfloor}\right]\right]\cr
&& \qquad = \mbb{E}\left[\e_\lambda(Z_k(\tau_k))\left[\mbb{E}\left[M_{k,\lambda}(\tau_k+\delta_k)\Big|\mcr{F}_{\lfloor\gamma_k\tau_k\rfloor}\right] - M_{k,\lambda}(\tau_k)\right]\right]\cr
&& \qquad = 0,
\eeqnn
which implies that $I_3 = 0.$ The result holds.
\qed

\bcorollary\label{c3.4}
Suppose that Condition \ref{condition} holds and $h \in C_b(\mbb{R}_+^2)^+$. Then for any $\lambda: = (\lambda_1, \lambda_2) \in \mbb{R}_+^2,$ we have
\beqnn
\lim_{k \rightarrow \infty} L^\lambda_{\tau_k, \delta_k}(Z_k) = 0,
\eeqnn
where $ L^\lambda_{\tau_k, \delta_k}(Z_k) := \mbb{E}\left[\left|e^{-\lambda_1X_k(\tau_k + \delta_k)} - e^{-\lambda_1X_k(\tau_k)}\right|^2 + \left|e^{-\lambda_2Y_k(\tau_k + \delta_k)} - e^{-\lambda_2Y_k(\tau_k)}\right|^2\right].$
\ecorollary

\proof
The result follows by taking $\lambda = (\lambda_1, 0)$ and $\lambda = (0, \lambda_2)$ in Proposition \ref{p1.3}.
\qed

Similar to the proof of \cite[Theorem 3.6]{L20}, we get the following result.

\bproposition\label{ptight}
Suppose that Condition \ref{condition} holds and $h \in C_b(\mbb{R}_+^2)^+$. Let $Z_k(0) = (X_k(0), Y_k(0))$ be the initial value satisfying  $\sup_k \mbb{E}[X_k(0) + Y_k(0)] < \infty$. Then the process $\{(Z_k(t))_{t \ge 0}\}_{k \ge 1} = \{(X_k(t), Y_k(t))_{t \ge 0}\}_{k \ge 1}$ is tight on ${\mbb{D}}([0, \infty), D).$
\eproposition

\proof
By Aldous's criterion, it suffices to show that, for any $\epsilon > 0,$
\beqlb\label{tight}
\lim_{k \rightarrow \infty}\mbb{P}\left[\left\|Z_k(\tau_k+\delta_k) - Z_k(\tau_k)\right\|_2 > \epsilon \right] = 0,
\eeqlb
where $\|\cdot\|_2$ is the $L^2$ norm on $D.$
For any $a :=(a_1, a_2),\ b:=(b_1, b_2) \in D$ satisfying $\|a - b\|_2 > \epsilon$,  we have $|a_1 - b_1| \wedge |a_2 - b_2| > \epsilon/2.$ Then for a fixed constant $M > 0,$ by taking $0 \le \|a\|_2, \|b\|_2 \le M,$ one sees that
\beqnn
|\e^{-\lambda_1a_1} - \e^{-\lambda_1b_1}|^2 + |\e^{-\lambda_2a_2} - \e^{-\lambda_2b_2}|^2 \ge \left(\frac{1}{2}(\lambda_1\wedge \lambda_2) \epsilon e^{-(\lambda_1 + \lambda_2) M}\right)^2.
\eeqnn
By Proposition \ref{p1.3}, it is easy to see that
\beqnn
\ar\ar\mbb{P}\left\{\|Z_k(\tau_k + \delta_k) - Z_k(\tau_k)\|_2 > \epsilon; \|Z_k(\tau_k)\|_2\vee\|Z_k(\tau_k + \delta_k)\|_2 \le M \right\}\cr
\ar\ar \quad \le \left(\frac{1}{2}(\lambda_1\wedge \lambda_2) \epsilon e^{-(\lambda_1 + \lambda_2) M}\right)^{-2} L^\lambda_{\tau_k, \delta_k}(Z_k) \rightarrow 0
\eeqnn
as $ k \rightarrow \infty.$ Further, by Proposition \ref{moment}, we have
\beqnn
\mbb{P}\left[\|Z_k(\tau_k + \delta_k)\|_2 \ge M\right] &\le& \mbb{P}\left[X_k(\tau_k + \delta_k) \ge \frac{M}{2}\right] + \mbb{P}\left[Y_k(\tau_k + \delta_k) \ge \frac{M}{2}\right]\cr
&\le& 2\frac{\sup_{0 \le t \le T}\mbb{E}[X_k(t) + Y_k(t)]}{M} \le \frac{K}{M}.
\eeqnn
Similarly, we get
\beqnn
\mbb{P}\left[\|Z_k(\tau_k)\|_2 \ge M\right] \le \frac{K}{M}.
\eeqnn
As a result,
\beqnn
\ar\ar\mbb{P}\left[\left\|Z_k(\tau_k+\delta_k) - Z_k(\tau_k)\right\|_2 > \epsilon \right] \cr
\ar\ar\qquad\le \mbb{P}\left[\|Z_k(\tau_k + \delta_k) - Z_k(\tau_k)\|_2 > \epsilon; \|Z_k(\tau_k)\|_2\vee\|Z_k(\tau_k + \delta_k)\|_2 \le M \right]\cr
\ar\ar\qquad\quad +\mbb{P}\left[\|Z_k(\tau_k + \delta_k) - Z_k(\tau_k)\|_2 > \epsilon; \|Z_k(\tau_k + \delta_k)\|_2 \ge M\right]\cr
\ar\ar\qquad\quad + \mbb{P}\left[\|Z_k(\tau_k + \delta_k) - Z_k(\tau_k)\|_2 > \epsilon; \|Z_k(\tau_k)\|_2 \ge M\right]\cr
\ar\ar\qquad\le \mbb{P}\left[\|Z_k(\tau_k + \delta_k) - Z_k(\tau_k)\|_2 > \epsilon; \|Z_k(\tau_k)\|_2\vee\|Z_k(\tau_k + \delta_k)\|_2 \le M \right]\cr
\ar\ar\qquad\quad +\mbb{P}\left[\|Z_k(\tau_k + \delta_k)\|_2 \ge M\right] + \mbb{P}\left[\|Z_k(\tau_k)\|_2 \ge M\right]
\eeqnn
goes to $0$ as $k \rightarrow \infty$ and $M \rightarrow \infty,$ which implies \eqref{tight}. The result follows.
\qed

% Similar to the proof of Lemma 2.2 of \cite{L20}, we get the following result.

\blemma\label{f^mn}
For any $f \in C_b^2(D)$, there exists a sequence of functions $f^{m,n} \in C_0^2(D)$ such that $f^{m,n} \rightarrow f$, $f^{m,n}_1 \rightarrow f_1$ and $f^{m,n}_{11} \rightarrow f_{11}$ uniformly on any bounded subset of $D$ as $m, n \rightarrow \infty$, where $f^{m,n}_1:=\frac{\partial f^{m,n}(x,y)}{\partial x}$, $f_1:=\frac{\partial f(x,y)}{\partial x}$, $f^{m,n}_{11}:=\frac{\partial^2 f^{m,n}(x,y)}{\partial x^2}$ and $f_{11}:=\frac{\partial^2 f(x,y)}{\partial x^2}$.
\elemma

\proof
For any nonnegative function $f \in C_b^2(D),$ we define
\beqnn
f^{m,n}(x, y) =
\begin{cases}
	f(x, y), &\quad (x, y) \in [0, m]\times[0,n]\cap D,\\
	f(x, y)\left[1 - 2\int_m^x\rho(2(z-m)-1)d z\right], &\quad (x, y) \in [m,m+1]\times[0,n]\cap D;\\
	0, &\quad \text{others,}
\end{cases}
\eeqnn
where $\rho$ is the mollifier defined by
\beqnn
\rho(x) = \Lambda\exp\{-1/(1-x^2)\}1_{\{|x| < 1\}}
\eeqnn
with $\Lambda$ being the constant such that $\int_{\mbb{R}}\rho(x) d x = 1.$ It is easy to see that $f^{m,n} \in C_0^2(D)$. Notice that, for $(x, y) \in [m,m+1]\times[0,n]\cap D,$
\beqnn
f^{m,n}_1(x, y) = f_1(x, y) - \frac{d}{d x}\left[2f(x, y)\int_m^x \rho(2(z-m)-1)d z \right]
\eeqnn
and
\beqnn
f^{m,n}_{11}(x, y) = f_{11}(x, y) - \frac{d^2}{d x^2}\left[2f(x, y)\int_m^x \rho(2(z-m)-1)d z \right].
\eeqnn
Let $D^b$ be a fixed bounded subset of $D.$ Then we have
\beqnn\sup_{(x, y) \in D^b}\left[|f^{m,n}(x, y) - f(x, y)| + \left|f^{m,n}_1(x, y) - f_1(x, y)\right| + \left|f^{m,n}_{11}(x, y) - f_{11}(x, y)\right|\right] \rightarrow 0
\eeqnn
as $m, n \rightarrow \infty.$ The result follows. \qed

Now we are ready to give the existence of the solution to \eqref{1.1a}--\eqref{1.1b} for the case of  $h \in C_b(D)^+$.

{\bf Proof of Theorem~\ref{existence}}\quad
Let $P^{(k)}$ be the distributions of $Z_k$ on ${\mbb{D}}([0, \infty), D)$.  By Proposition \ref{ptight}, the sequence of processes $\{Z_k\}_{k \ge 1}$ is relatively compact. Then there are a probability measure $Q$ and a subsequence $P^{(k_i)}$  on ${\mbb{D}}([0, \infty), D)$ such that $Q = \lim_{i \rightarrow \infty}P^{(k_i)}.$ By Skorokhod representative theorem, there exists a probability space $(\tilde{\Omega}, \tilde{\mcr{F}}, \tilde{\mbb{P}})$ on which are defined c\`{a}dl\`{a}g processes $(\tilde{Z}(t))_{t \ge 0}$ and $(\tilde{Z}_{k_i}(t))_{t \ge 0}$ such that the distribution of $\tilde{Z}$ and $\tilde{Z}_{k_i}$ on ${\mbb{D}}([0, \infty), D)$ are $Q$ and $P^{(k_i)},$ respectively, and  $\lim_{i \rightarrow \infty} \tilde{Z}_{k_i} = \tilde{Z},\ \tilde{\mbb{P}}$-almost surely.

Now it suffices to show that $(\tilde{Z}(t))_{t \ge 0}$ satisfies the following martingale problem: for any $f \in C_b^2(D),$ we have
\beqlb\label{MP}
f(\tilde{Z}(t)) = f(\tilde{Z}(0)) + \int_0^t Af(\tilde{Z}(s)) d s + \text{local\ mart.}
\eeqlb
Let $f(z) = \e_\lambda(z)$ for any $z \in D.$ By \eqref{eZk}, we get
\beqnn
\e_\lambda(\tilde{Z}_{k_i}(t)) = \e_\lambda(\tilde{Z}_{k_i}(0)) + \int_0^{\lfloor\gamma_{k_i}t\rfloor/\gamma_{k_i}} A_{k_i}\e_\lambda(\tilde{Z}_{k_i}(s))d s + M_{{k_i},\lambda}(t).
\eeqnn
One sees that
\beqnn
\ar\ar\int_0^{\lfloor\gamma_{k_i}t\rfloor/\gamma_{k_i}}\left|A_{k_i} \e_\lambda(\tilde{Z}_{k_i}(s)) - A\e_\lambda(\tilde{Z}(s))\right|d s\cr
\ar\ar\quad\le \int_0^{\lfloor\gamma_{k_i}t\rfloor/\gamma_{k_i}}\left|A_{k_i} \e_\lambda(\tilde{Z}_{k_i}(s)) - A\e_\lambda(\tilde{Z}_{k_i}(s))\right|d s\cr
\ar\ar\qquad\quad + \int_0^{\lfloor\gamma_{k_i}t\rfloor/\gamma_{k_i}}\left|A\e_\lambda(\tilde{Z}_{k_i}(s)) - A\e_\lambda(\tilde{Z}(s))\right|d s =: I_{k_i}^1 + I_{k_i}^2.
\eeqnn
Then $I_{k_i}^1 \rightarrow 0$ as $i \rightarrow \infty$ by Theorem \ref{elimit}. On the other hand, let
$$
C_X := \{t > 0: \tilde{P}(\tilde{Z}(t-) = \tilde{Z}(t)) = 1\}.
$$
Then the set $\mbb{R}_+\backslash C_X$ is at most countable. Then we have $I_{k_i}^2 \rightarrow 0$ as $i \rightarrow \infty.$ Consequently, the process $(\tilde{Z}(t))_{t \ge 0}$ satisfies the martingale problem \eqref{MP} when $f(z) =\e_\lambda(z).$

Let $f \in C_0^2(D)$ be fixed, and $E_0$ be the linear hull of $\{\e_\lambda(z): \lambda \in \mbb{R}_+^2\}.$ By Stone-Weierstrass Theorem and \eqref{generator}, there exists a sequence of functions $f_n \in E_0$ such that $Af_n(z) \rightarrow Af(z)$ uniformly on each bounded subset of $D$ as $n \rightarrow \infty$. As a linear span of $\{\e_\lambda(z)\},$ we have
\beqlb\label{fnZ}
f_n(\tilde{Z}(t)) = f_n(\tilde{Z}(0)) + \int_0^t Af_n(\tilde{Z}(s)) d s + \text{local\ mart.}
\eeqlb
Let $\tilde{\tau}_N := \inf\{t > 0: \tilde{X}(t) \ge N\ \text{or}\ \tilde{Y}(t) \ge N\}.$ Then $\tilde{\tau}_N \rightarrow \infty$ almost surely as $N \rightarrow \infty$ by Proposition \ref{moment} and Fatou's lemma. Replacing $t$ with $t \wedge \tilde{\tau}_N,$ and taking limits as $n \rightarrow \infty$ on both sides of \eqref{fnZ}, we then have
\beqlb\label{MPf}
f(\tilde{Z}(t \wedge \tilde{\tau}_N)) = f(\tilde{Z}(0)) + \int_0^t Af(\tilde{Z}(s \wedge \tilde{\tau}_N-)) d s + \text{mart.}
\eeqlb

Next, for the general function $f \in C_b^2(D),$ by Lemma \ref{f^mn}, there exists a sequence functions $f^{m,n} \in C_0^2(D)$ such that $Af^{m,n}(z) \rightarrow Af(z)$ uniformly on each bounded subset of $D$ as $m, n \rightarrow \infty.$ Similar to the above, \eqref{MPf} holds for any $f \in C_b^2(D).$
Letting $N \rightarrow \infty,$ one can see that $(\tilde{Z}(t))_{t \ge 0}$ satisfies the martingale problem \eqref{MP}, which implies that $(\tilde{Z}(t))_{t \ge 0}$ is a weak solution to \eqref{1.1a}--\eqref{1.1b}. The result follows.
\qed

\begin{lemma}\label{l210}
	Assume that $h \in C_b(\mbb{R}_+^2)^+.$ Let $(X(t), Y(t))_{t \ge 0}$ be the solution to \eqref{1.1a}-\eqref{1.1b} with $\mbb{E}[X(0) + Y(0)] < \infty$. Then for any $T > 0$ we have
	\beqnn
	\sup_{0 \le t \le T}\mbb{E}[X(t) + Y(t)] < \infty.
	\eeqnn
\end{lemma}

\proof
Taking $f(x, y) = x + y$ and $Z(t) = (X(t), Y(t)),$ by \eqref{MPf} we have
\beqnn
\mbb{E}[X(t\wedge \tilde{\tau}_N) + Y(t\wedge \tilde{\tau}_N)] = \mbb{E}[X(0)+Y(0)] + \int_0^t \mbb{E}\left[Af(Z(s \wedge \tilde{\tau}_{N-})) \right]\d s.
\eeqnn
By \eqref{generator} and $h \in C_b(\mbb{R}_+^2)^+,$ one sees that
\beqnn
\mbb{E}[X(t\wedge \tilde{\tau}_N) + Y(t\wedge \tilde{\tau}_N)] \ar\le\ar \mbb{E}[X(0)+Y(0)]\cr
\ar\ar \!+ \int_0^t \mbb{E}\left[-b X(s\!\wedge\! \tilde{\tau}_N)+ \left(\!\sup_{(x,y) \in \mbb{R}_+^2}\!h(x,y)\! \int_{\mbb{N}^{-1}}\xi n(\d \xi)\!\right)\! Y(s\!\wedge\! \tilde{\tau}_N)\right]\! \d s \cr
\ar\le\ar \mbb{E}[X(0)+Y(0)] + K\int_0^t \mbb{E}[X(s\wedge \tilde{\tau}_N)+Y(s\wedge \tilde{\tau}_N)] \d s.
\eeqnn
Then the result follows by Gronwall's inequality and letting $N \rightarrow \infty.$
\qed

{\bf Proof of Theorem \ref{t1.1}}\quad
By \cite[Theorems 5.1 and 5.2]{DL06} and \cite[Corollary 5.2]{FL10}, there is a unique positive strong solution to \eqref{1.1a}, which is a CB-process. Let $(X(t))_{t \ge 0}$ be the unique positive strong solution to \eqref{1.1a}. The pathwise uniqueness of solution to \eqref{1.1b} can be constructed by path stitching method. Let $\kappa_0 = 0.$ Then $Y(\kappa_0) = Y(0).$ Given $\kappa_{k - 1} \ge 0$, $Y(\kappa_{k - 1})\ge 0$ and the process $(X(t))_{t \ge 0},$ we first define
\beqnn
\kappa_k = \kappa_{k - 1} + \inf\left\{t > 0: \int_{\kappa_{k - 1}}^{t + \kappa_{k - 1}}\int_0^{Y(\kappa_{k - 1})}\int_0^{h(X(s-), Y(\kappa_{k - 1}))}\int_{\mbb{N}^{-1}}N(\d s, \d u, \d r, \d \xi) = 1\right\}.
\eeqnn
Then we define $Y(t) = Y(0)$ for any $t \in [\kappa_0, \kappa_1)$ and
\beqnn
\Delta_k :=  \int_{\kappa_{k - 1}}^{\kappa_{k}}\int_0^{Y(\kappa_{k - 1})}\int_0^{h(X(s-), Y(\kappa_{k - 1}))}\int_{\mbb{N}^{-1}}\xi N(\d s, \d u, \d r, \d \xi).
\eeqnn
For any $k = 1, 2, \cdots,$ let $Y(t) := Y(\kappa_{k - 1}) + \Delta_k$ for any $t \in [\kappa_k, \kappa_{k + 1}),$ which uniquely determine the behavior of the trajectory $t \rightarrow Y(t)$ on the time interval $[\kappa_k, \kappa_{k + 1}),\ k = 1, 2, \cdots.$ Let $\kappa := \lim_{k \rightarrow \infty}\kappa_k.$ Then $(Y(t))_{t \ge 0}$ is the pathwise-unique solution to \eqref{1.1b} up to time $\kappa.$ Notice that $\kappa_k$ is the time of $k$-th jump of $Y.$ By Lemma \ref{l210}, one sees that
\beqnn
\mbb{E}\left[\int_0^{t\wedge \kappa_k}\int_0^{Y(s-)}\int_0^{h(X(s-), Y(s-))}\int_{\mbb{N}^{-1}}N(\d s, \d u, \d r, \d \xi)\right] \le Kt \sup_{0 \le s \le t}\mbb{E}[Y(s)] < \infty.
\eeqnn
It follows that $\mbb{P}(\kappa > t) = 1$ for any $t \ge 0,$ which implies that $\mbb{P}(\kappa = \infty) = 1.$  Then $(X(t), Y(t))_{t \ge 0}$ is the pathwise-unique solution to \eqref{1.1a}-\eqref{1.1b}. The result follows.
\qed

{\bf Proof of Theorem \ref{texistence}}\quad
Let $h_m(x, y) := h(x\wedge m, y\wedge m).$ Then $h_m$ is bounded for any $m \ge 1$ and $h_m \rightarrow h$ as $m \rightarrow \infty$. By Theorems \ref{existence} and \ref{t1.1}, there exists a unique strong solution $(X_m(t), Y_m(t))_{t \ge 0}$ to the following stochastic integral equation system:
\beqlb\label{m}
\begin{cases}
	X(t)=X(0)-b\int_0^t X(s)\,\d s+\int_0^t\sqrt{2cX(s)}\,\d B(s)+\int_0^t\int_0^{X(s-)}\int_0^\infty \xi\,\tilde{M}(\d s,\d u,\d \xi),\\
	Y(t)=Y(0)+\int_0^t\int_0^{Y(s-)}\int_0^{h_m(X(s-),Y(s-))}\int_{\mbb{N}^{-1}} \xi\,N(\d s,\d u,\d r,\d \xi).
\end{cases}
\eeqlb
In fact, $(X_m(t))_{t \ge 0}$ is the unique strong solution to \eqref{1.1a} independent with $m$, which is written as $(X(t))_{t \ge 0}$ in the following. Let $\tau_m^X :=\inf\{t > 0: X(t)\ge m\},\ \tau_m^Y :=\inf\{t > 0: Y_m(t) \ge m\}$ and $\tau_m = \tau_m^X \wedge \tau_m^Y.$ Then $0 \le X(t) < m$ and $0 \le Y_m(t) < m$ for $0 \le t < \tau_m,$ and $(X(t), Y_m(t))$ satisfies \eqref{1.1a}-\eqref{1.1b} for $0 \le t < \tau_m.$
For $n \ge m \ge 1$, let
\beqnn
Y_m(\tau_m) = Y_m(\tau_m-) + \int_{\{\tau_m\}}\int_0^{Y_m(\tau_m-)}\int_0^{h_n(X(\tau_m-),Y_m(\tau_m-))}\int_{\mbb{N}^{-1}} \xi N(\d s,\d u,\d r,\d \xi).
\eeqnn
There exists a unique strong solution $(X(t), \tilde{Y}(t))_{t \ge \tau_m}$ to \eqref{1.1a} and
\beqnn
Y(t)=Y_m(\tau_m)+\int_{\tau_m}^t\int_0^{Y(s-)}\int_{\tau_m}^{h_n(X(s-),Y(s-))}\int_{\mbb{N}^{-1}} \xi\,N(\d s,\d u,\d r,\d \xi).
\eeqnn
Let $Y'(t) = Y_m(t)$ for $0 \le t < \tau_m$ and $Y'(t)= \tilde{Y}(t)$ for $t \ge \tau_m.$ Then it is a solution to \eqref{m} by changing $m$ to $n.$ By the strong uniqueness we get $(X(t), Y'(t))_{t \ge 0} = (X(t), Y_n(t))_{t \ge 0}$ almost surely. In particular, we infer $Y_n(t) = Y_m(t) < m$ for $0 \le t < \tau_m$. Consequently, the sequence $\{\tau_m\}$ is non-decreasing.
On the other hand, by \eqref{m} it is easy to check that $\mbb{E}[X(t\wedge \tau_m^X)] \le \mbb{E}[X(0)]e^{Kt},$ where $K$ is a constant independent with $m.$ Then we have
$\tau_m^X \rightarrow \infty$ almost surely as $m \rightarrow \infty$. Let $\tau = \lim_{m \rightarrow \infty} \tau_m = \lim_{m \rightarrow \infty}\tau_m^Y$. Let $Y(t) = Y_m(t)$ for all $0 \le t < \tau_m$ and $m \ge 1$. It is easily seen that $(X(t), Y(t))_{t \in [0, \tau)}$ is a unique strong solution to \eqref{1.1a}--\eqref{1.1b} up to $\tau$. For $t \ge \tau,$ let $(X(t), Y(t)) = (X(t), \infty).$ The result follows.
\qed

\begin{comment}
\btheorem\label{comparison}
 Suppose that $(X_1(t), Y_1(t))_{t \ge 0}$ and $(X_2(t), Y_2(t))_{t \ge 0}$ are two positive solutions to \eqref{1.1a}-\eqref{1.1b} with $\mbb{P}(X_1(0) \le X_2(0),\ Y_1(0) \le Y_2(0)) = 1.$ Let $h \in C(\mbb{R}_+^2)^+$. Then we have $\mbb{P}(X_1(t) \le X_2(t),\ Y_1(t) \le Y_2(t)\ \text{for\ all}\ t \ge 0) = 1.$
\etheorem
\proof
 By \cite[Theorem 8.4]{L20}, we have $\mbb{P}(X_1(t) \le X_2(t)\ for\ all\ t \ge 0) = 1.$ For $n \ge 0$ let $\phi_n$ be the function defined as in the proof of Theorem \ref{t1.1}. Let $\psi_n(z) = \phi_n(z \vee 0)$ for $z \in \mbb{R}.$ Then $\psi_n(z) \rightarrow z_+ := z \vee 0$ increasingly as $n \rightarrow \infty.$  For any $k \in \mbb{N}_+,$ we define $\kappa_k^1 = \inf\{t \ge 0: X_1(t) + Y_1(t) \ge k\},$ $\kappa_k^2 = \inf\{t \ge 0: X_2(t) + Y_2(t) \ge k\}$ and $\kappa_k = \kappa_k^1 \wedge \kappa_k^2.$ Let $\kappa := \lim_{k \rightarrow \infty}\kappa_k.$ Then $Y_1(t) = Y_2(t) = \infty$ almost surely for any $t \ge \kappa.$ Let $\Delta(t) = Y_1(t) - Y_2(t)$ for $t \in [0, \kappa).$ Similar to the proof of Theorem \ref{t1.1}, we get $\mbb{E}[(\Delta(t\wedge \kappa_k))_+] = 0$ for every $t \ge 0$ by Gronwall's inequality. Then by taking $k \rightarrow \infty$ and using Fatou's lemma we see that $\mbb{E}[(Z(t))_+] = 0$ for every $t \in [0, \kappa),$ then $\mbb{P}(Y_1(t) \le Y_2(t) \ \text{for\ all}\ t \ge 0) = 1$ by the right continuity of the processes.
\qed
\end{comment}

\section{Large time behaviors}

\subsection{Foster-Lyapunov criteria for extinction}

In this subsection, we mainly discuss the extinction behavior of such processes under $b\geq0$.  Define $\tau_0 = \inf\{t > 0: X(t) = 0\ \text{and}\ Y(t) = 0\}.$ Moreover, we separately define the extinction time of $X, Y$ as  $\tau_0^X := \inf\{t > 0: X(t) = 0\}$ and $\tau_0^Y := \{t>0: Y(t) = 0\}.$ Then we have $\tau_0 = \tau_0^X \vee \tau_0^Y.$ For the extinction behavior of the process $X$, we introduce the so called Grey's condition:
\bcondition\label{grey cond}
There is some constant $\theta>0$ so that $\phi_1(z)>0$ for $z\geq\theta$ and $\int_\theta^{\infty}\phi_1^{-1}(z)\d z<\infty$, where $\phi_1$ is given by \eqref{phi}.
\econdition

%Define $\tau^{X,-}_x:=\inf\{t>0: X_t<x\}$ and $\tau^{X,+}_x=\inf\{t>0: X_t>x\}$ for $x\geq0$.

Since $b\geq0$, under Condition \ref{grey cond}, one can see that $\mbb{P}_x(\tau^{X}_0<\infty)=1$ for all $x>0$; see, e.g., \cite[Corollary 3.8]{L20}. In the following, we present a Foster-Lyapunov criteria-type result for the process $(X,Y)$.
For $z_1=(x_1,y_2),z_2=(x_2,y_2)\in D$, we say $z_1\succeq z_2$ if $x_1\geq x_2$ and $y_1\geq y_2$. Let $\tilde{z} := (\tilde{x}, \tilde{y}) \succeq z_0$ and $(X(t), Y(t))_{t \ge 0}$ be the mixed state branching process satisfying \eqref{1.1a}--\eqref{1.1b} with initial value $z_0.$ We define stopping time $\sigma_{\tilde{z}} = \inf\{t > 0: X(t) \ge \tilde{x} \ \text{or}\ Y(t)\ge \tilde{y}\}.$ It is easy to see that $X_{t \wedge \sigma_{\tilde{z}}-} \le \tilde{x}$ and $Y_{t \wedge \sigma_{\tilde{z}}-} \le \tilde{y}.$

\begin{theorem}\label{foster-lyapunov}
	Let $(Z(t))_{t\geq0}=(X(t), Y(t))_{t \ge 0}$ be the mixed state branching process satisfying \eqref{1.1a}--\eqref{1.1b} with initial value $z_0 =(x_0,y_0)\in D.$ Suppose that $\phi_1(\lambda_1) > 0$ and $\phi_2(\lambda_2) > 0$ for any $\lambda :=(\lambda_1, \lambda_2) \in (0, \infty)^2.$
	Then we have $\mbb{P}_{z_0}\{\tau_{0}<\infty\} = 1.$
\end{theorem}
\proof
It suffices to prove the case of $z_0 \in D\backslash (0, 0).$ The proof is inspired by \cite[Lemma 4.1]{LYZ19}.
By It\^{o}'s formula, we have
\beqlb
\mrm{e}_\lambda(Z(t\wedge \tau_0\wedge \sigma_{\tilde{z}})) = \mrm{e}_\lambda(z_0) + \int_0^{t\wedge \tau_0\wedge \sigma_{\tilde{z}}} A\mrm{e}_\lambda(Z(s-)) \d s + mart.
\eeqlb
Taking expectations on both sides, we have
\beqnn
\mbb{E}_{z_0}\left[\mrm{e}_\lambda(Z(t\wedge \tau_0\wedge \sigma_{\tilde{z}}))\right] = \mrm{e}_\lambda(z_0) + \int_0^t \mbb{E}_{z_0} \left[ A\mrm{e}_\lambda(Z(s-)) 1_{\{s< \tau_0\wedge \sigma_{\tilde{z}}\}}\right]\d s,
\eeqnn
which implies that
\beqnn
\d(\mbb{E}_{z_0}\left[\mrm{e}_\lambda(Z(t\wedge \tau_0\wedge \sigma_{\tilde{z}}))\right]) = \mbb{E}_{z_0} \left[ A\mrm{e}_\lambda(Z(t-)) 1_{\{t< \tau_0\wedge \sigma_{\tilde{z}}\}}\right]\d t.
\eeqnn
Recall that $\phi_1(\lambda_1) > 0$ and $\phi_2(\lambda_2) > 0$ for all $\lambda \in (0, \infty)^2.$ Then for all $\tilde{z}=(\tilde{x},\tilde{y})\in D$ with $\tilde{z}\succeq z_0$ and $\lambda \in (0, \infty)^2$, there exists a constant $d_{z_0, \tilde{z}, \lambda}>0$ such that for all $z=(x,y)\in D$ with $z_0\preceq z\preceq\tilde{z}$,
\beqlb\label{Lyapnuov cond}
x\phi_1(\lambda_1)+h(x,y)y\phi_2(\lambda_2)\geq d_{z_0, \tilde{z},\lambda}.
\eeqlb	
Then by integration by parts,
\beqnn
\ar\ar\int_0^\infty\e^{-d_{z_0, \tilde{z}, \lambda}t}\mbb{E}_{z_0}\left[ A\mrm{e}_\lambda(Z(t)) 1_{\{t< \tau_0\wedge \sigma_{\tilde{z}}\}}\right]\d t\cr
\ar\ar\qquad = \int_0^\infty\e^{-d_{z_0, \tilde{z}, \lambda}t}\d (\mbb{E}_{z_0}\left[\mrm{e}_\lambda(Z(t\wedge \tau_0\wedge \sigma_{\tilde{z}}))\right])\cr
\ar\ar\qquad = d_{z_0, \tilde{z}, \lambda}\int_0^\infty\e^{-d_{z_0, \tilde{z}, \lambda}t}  \mbb{E}_{z_0}\left[\mrm{e}_\lambda(Z(t\wedge \tau_0\wedge \sigma_{\tilde{z}}))\right]\d t - \mrm{e}_\lambda(z_0).
\eeqnn
Moreover, by \eqref{Lez} and \eqref{Lyapnuov cond} we have
\beqnn
\ar\ar\int_0^\infty \e^{-d_{z_0, \tilde{z}, \lambda}t}\mbb{E}_{z_0}\left[ A\mrm{e}_\lambda(Z(t)) 1_{\{t< \tau_0\wedge \sigma_{\tilde{z}}\}}\right]\d t\cr
\ar\ar\qquad \ge d_{z_0, \tilde{z}, \lambda} \int_0^\infty \e^{-d_{z_0, \tilde{z}, \lambda}t} \mbb{E}_{z_0}\left[\mrm{e}_\lambda(Z(t)) 1_{\{t< \tau_0\wedge \sigma_{\tilde{z}}\}}\right]\d t.
\eeqnn
It follows that
\beqnn
\e_\lambda(z_0) &\le& d_{z_0, \tilde{z}, \lambda} \int_0^\infty e^{-d_{z_0, \tilde{z}, \lambda}t} \mbb{E}_{z_0}\left[\mrm{e}_\lambda(Z(\tau_0\wedge \sigma_{\tilde{z}})) 1_{\{t\ge \tau_0\wedge \sigma_{\tilde{z}}\}}\right]\d t\cr
&\le& \mbb{P}_{z_0}\{\tau_0\le \sigma_{\tilde{z}}\} + \sup_{z \succeq \tilde{z}}[\e^{- \lambda_1 x} + \e^{- \lambda_2 y}]\cr
&\le& \mbb{P}_{z_0}\{\tau_0< \infty\} + [\e^{- \lambda_1 \tilde{x}} + \e^{- \lambda_2 \tilde{y}}].
\eeqnn
Taking $\tilde{x}, \tilde{y} \rightarrow \infty,$ we get $\mbb{P}_{z_0}\{\tau_0< \infty\} \ge \mrm{e}_\lambda(z_0),$ which holds for any  $\lambda\in(0,\infty)^2$.
The result follows by letting $\lambda\rightarrow(0,0).$
\qed
\begin{remark}
	The processes $(X(t))_{t \ge 0}$ and
	$(Y(t))_{t \ge 0}$ are independent when $h$ is a positive constant. In this case, one can check that $\mbb{P}(\tau_0^X < \infty) = 1$ when $\phi_1(\lambda_1) > 0$ for any $\lambda_1 > 0$, and $\mbb{P}(\tau_0^Y < \infty) = 1$ if $\phi_2(\lambda_2) > 0$ for any $\lambda_2 > 0.$
\end{remark}
\bcorollary\label{C4.1}
Assume that $b \ge 0$, $R_1 := \int_{\mbb{N}^{-1}}\xi n(\d \xi) < 0$ and Condition \ref{grey cond} holds. Then we have $\mbb{P}_{z_0}\{\tau_{0}<\infty\} = 1.$
\ecorollary

\proof
By Condition \ref{grey cond} and $b \ge 0,$ one sees that $\phi_1(\lambda_1) > 0$ for any $\lambda_1 > 0.$ Moreover, by the inequality $1 - e^{-\lambda_2 \xi} \le \lambda_2 \xi,$ we have $\phi_2(\lambda_2) \ge - R_1\lambda_2 > 0.$ The result follows by Theorem \ref{foster-lyapunov}.
\qed

\subsection{Exponential ergodicity in the $L^1$-Wasserstein distance}

The coupling method  is a powerful tool in the study of ergodicity of Markov processes.  We refer the reader to \cite{C05,L92} for the systematical study on this topic.

We say a couple process $(Z_1(t),Z_2(t))_{t\geq0}$ is called a {\it coupling} of $(Z(t))_{t\geq0}$ with transition semigroup $(P_t)_{t\geq0}$ if both $(Z_1(t))_{t\geq0}$ and $(Z_2(t))_{t\geq0}$ are Markov processes with transition semigroup $(P_t)_{t\geq0}$ (possibly with different initial distributions) and
$Z_1(t+\tau) = Z_2(t+\tau)$ for every $t\geq0$,
where
$$
\tau=\inf\{t\geq0:Z_1(t)=Z_2(t)\}.
$$
In this case, $(Z_1(t))_{t\geq0}$ and $(Z_2(t))_{t\geq0}$ are called the {\it marginal processes} of the coupling. Let $A$ and $\tilde{A}$ be infinitesimal generators of $(Z(t))_{t\geq0}$ and $(Z_1(t),Z_2(t))_{t\geq0}$, respectively.  Then $\tilde{A}$ is called the {\it coupling generator} and satisfies the marginal property, i.e., for any
$f,g\in{\cal{D}}(A)$,
\beqlb\label{defin coupling gener}
\tilde{A}u(z,\tilde{z})=Af(z)+Ag(\tilde{z})
\eeqlb
with $u(z,\tilde{z})=f(z)+g(\tilde{z})$.

A coupling $(Z_1(t),Z_2(t))_{t\geq0}$ is called {\it successful} if $\tau<\infty$ almost surely. A Markov process is said to have a {\it coupling property} if, for any initial distributions $\mu_1$ and $\mu_2$, there exists a successful coupling with marginal processes starting from $\mu_1$ and $\mu_2$, respectively. For any initial distribution $\mu$, let $\mathbb{P}^\mu$ be the distribution of this process with initial distribution $\mu$, and let $\mu P_t$  be  the marginal distribution of $ \mathbb{P}^\mu$ . It is known, see \cite{CG95,L92}, that the coupling property of the
process is equivalent to the following statement that
\beqnn
\mbox{For any initial distributions~~} \mu_1,  \mu_2,\  \lim_{t\rightarrow\infty}
\|\mu_1P_t-\mu_2P_t\|_{\mathrm{var}}=0,
 \eeqnn
where $\|\cdot\|_{\mathrm{var}}$ is the  total variational norm in the sense of $
\|\mu-\nu\|_{\mathrm{var}}:=\sup\{|\mu(A)-\nu(A)|: A: \text{Borel~set}\}.$

The total variational norm is a special case of {\it Wasserstein distances}. By ${\cal P}(D)$ we denote the space of all Borel probability measures over $D$. Given $\mu,\nu\in{\cal P}(D)$, a coupling $H$ of $(\mu,\nu)$ is a Borel probability measure on $D\times D$ which has marginals $\mu$ and $\nu$, respectively. We write ${\cal H}(\mu,\nu)$ for the collection of all such couplings. Let $d$ be a metric on $D$ such that $(D,d)$ is a complete separable metric space and define
$$
{\cal P}_d(D)=\Big\{\rho\in{\cal P}(D):\int_D d(z,0)\,\rho(\d z)<\infty\Big\}.
$$
Then the Wasserstein distance on ${\cal P}_d(D)$ is defined by
$$
W_d(\mu,\nu)=\inf\Big\{\int_{D\times D}d(z,\tilde{z})\,H(\d z,\d \tilde{z}): H\in{\cal H}(\mu,\nu)\Big\}.
$$
Moreover, it can be shown that this infimum is attained; see, e.g., \cite[Theorem 6.16]{V09}. More precisely, there exists $\tilde{H}\in{\cal{H}}(\mu,\nu)$ such that
$$
W_d(\mu,\nu)=\int_{D\times D}d(z,\tilde{z})\,H(\d z,\d \tilde{z}).
$$
 Taking $d(z,\tilde{z}) = 1_{\{z\neq \tilde{z}\}}$, then we have ${\cal P}_d(D)={\cal P}(D)$ and
$W_d(\mu,\nu)=\|\mu-\nu\|_{\mathrm{var}}$; see, e.g., \cite[pp.18]{C05}. By taking $d(z,\tilde{z})=\|z-\tilde{z}\|_1,$ ${\cal P}_d(D):={\cal P}_1(D)$ defined by
$$
{\cal P}_1(D):=\Big\{\rho\in{\cal P}(D):\int_D \|z\|_1\,\rho(\d z)<\infty\Big\},
$$
where $\|\cdot\|_1$ is the $L^1$ norm on $D.$
We say the corresponding Wasserstein distance $W_1$ is the {\it $L^1$-Wasserstein distance.}

\begin{definition}
	We say $(Z(t))_{t\geq0}$ on $D$ or its transition semigroup $(P_t)_{t\geq0}$ is {\it exponential ergodic in the $L^1$-Wasserstein distance} with rate $\lambda_0>0$ if its possesses a unique stationary distribution $\mu$ and there is a nonnegative function $\nu\mapsto C(\nu)$ on ${\cal P}(D)$ such that
	\beqlb\label{exponential ergodic}
	W_1(\nu P_t,\mu)\leq C(\nu)\e^{-\lambda_0 t},\quad t\geq0, \nu\in{\cal P}(D).
	\eeqlb
\end{definition}
By standard arguments, (\ref{exponential ergodic}) follows if ($P_t(z,\cdot):=\delta_zP_t$)
\beqlb\label{exponential contract}
W_1(P_t(z,\cdot),P_t(\tilde{z},\cdot))\leq C_0(z,\tilde{z})\e^{-\lambda_0t},\quad t\geq0
\eeqlb
for $C_0(z,\tilde{z})>0$ depending on $z$ and $\tilde{z}$; see, e.g., the proof of \cite[Theorem 3.5]{FJKR23}.

In the literature of the exponential ergodicity of branching processes, \cite{LM15,L20,FJKR23,CL21} obtained the results by making full use of the branching property. If the property fails, it has been shown that coupling methods are very effective; see, e.g., \cite{LLWZ23+,LW20}. In this subsection, we mainly consider that the rate function $h(\cdot,\cdot)$ satisfies for any $z=(x,y)\in D$ that
\beqlb\label{function h}
h(x,y)=r+xm(y)
\eeqlb
with $r>0$ and $m(\cdot)\in C(\mbb{R}_+)$ is nonnegative.

This subsection consists two parts. First, we need to construct a coupling for the mixed state branching process with interactions $(Z(t))_{t\geq0}=(X(t),Y(t))_{t\geq0}$. Second, we construct a proper function $F(x,y,\tilde{x},\tilde{y})=F(|x-\tilde{x}|+|y-\tilde{y}|)$ for any $z=(x,y),\tilde{z}=(\tilde{x},\tilde{z})\in D$ satisfying two properties:

(i) the exponential contraction property: for any $z=(x,y),\tilde{z}=(\tilde{x},\tilde{y})\in D$,
$$
\tilde{A} F(|x-\tilde{x}|+|y-\tilde{y}|)\leq-\lambda F(|x-\tilde{x}|+|y-\tilde{y}|)
$$
with some $\lambda>0$, where $\tilde{A}$ denotes the coupling generator.

(ii) control the $L^1$-Wasserstein distance in the
sense that
$$
F(|x-\tilde{x}|+|y-\tilde{y}|)\asymp
|x-\tilde{x}|+|y-\tilde{y}|,
$$
here, $f\asymp g$ means that there is a constant $C_1\geq1$ such that,
$$
C_1^{-1}f(\cdot)\leq g(\cdot)\leq C_1f(\cdot).
$$
Then we can deduce (\ref{exponential ergodic}) from (i)--(ii).

Recalling that the generator $A$ of $(Z(t))_{t\geq0}=(X(t),Y(t))_{t\geq0}$ is given by \eqref{generator} for any $f\in C_b^{2,1}(\mbb{R}_+^2)$. Let ${\cal D}(A)$ denote the linear space consisting of functions $f\in C_b^{2,1}(\mbb{R}_+^2)$ such that the two integrals on the right-hand side of (\ref{generator}) are convergent and define continuous functions on $D$. To study the coupling and ergodicity of the process, we begin with the construction of a new coupling operator. we will combine the coupling by reflection for the Brownian motion and the synchronous coupling for Poisson random measures. Here the coupling by reflection for
Brownian motion means that we will take $(-B(t))_{t\geq0}$ (which is regarded as a reflection of $(B(t))_{t\geq0}$) for the process $(X(t))_{t\geq0}$ before two marginal processes meet. To explain the meaning of the synchronous coupling for Poisson random measures, we will use the viewpoint from the coupling operator. Set $z=(x,y),\ \tilde{z}=(\tilde{x},\tilde{y})\in D$ with $x\geq\tilde{x}$. Recall that $\gamma(x,y)=yh(x,y)$. The jump system corresponding to the synchronous coupling for the operator $A$ is given by
\beqnn
(x,\tilde{x}) \rightarrow \left\{
\begin{aligned}
	&	(x+\xi,\tilde{x}+\xi), ~~\tilde{x}m(\d\xi),\\
	&	(x+\xi,\tilde{x}),~~~~~~~
	(x -\tilde{x})m(\d\xi)
\end{aligned}
\right.
\eeqnn
and
\beqnn
(y,\tilde{y}) \rightarrow \left\{
\begin{aligned}
	&	(y+\xi,\tilde{y}+\xi), \qquad [\gamma(x,y)\wedge\gamma(\tilde{x},\tilde{y})]n(\d\xi),\\
	&	(y+\xi,\tilde{y}), \qquad ~~~~~
	[\gamma(x,y)-\gamma(\tilde{x},\tilde{y})]^+n(\d\xi),\\
	&	(y,\tilde{y}+\xi), \qquad ~~~~~
	[\gamma(x,y)-\gamma(\tilde{x},\tilde{y})]^-n(\d\xi),
\end{aligned}
\right.
\eeqnn
where $f^+$ ($f^-$) denotes the positive (negative) part of function $f$. Similarly, we can define the case that $x<\tilde{x}.$ We refer to \cite{C05} for the details of such coupling and other couplings for jump systems. We also refer to \cite[remark 2.4]{LW20} for similar discussions in the setting of nonlinear continuous state branching processes. Let $C^2(D\times D)$ be the set of continuous function $(x, y, \tilde{x}, \tilde{y})\mapsto F(x, y, \tilde{x}, \tilde{y})$ on $D\times D$ with continuous derivatives uo to $2$nd order on $x$ and $\tilde{x}.$
With the idea above in mind, we then define for any $F\in C^{2}(D\times D)$ and $x\geq\tilde{x}$ that
\begin{equation}\label{coupling gener}
\begin{split}
\tilde{A}F(x,y,\tilde{x},\tilde{y})
&= -bxF'_x-b\tilde{x}F'_{\tilde{x}}
+cxF''_{xx}+c\tilde{x}F''_{\tilde{x}\tilde{x}}-2c\sqrt{x\tilde{x}}F''_{x\tilde{x}}\cr
&~+\tilde{x}\int_0^\infty
[F(x+\xi,y,\tilde{x}+\xi,\tilde{y})-F(x,y,\tilde{x},\tilde{y})-\xi(F'_x+F'_{\tilde{x}})]
\,m(\d\xi)\cr
&~+(x-\tilde{x})\int_0^\infty
[F(x+\xi,y,\tilde{x},\tilde{y})-F(x,y,\tilde{x},\tilde{y})-\xi F'_x]
\,m(\d\xi)\cr
&~+[\gamma(x,y)\wedge\gamma(\tilde{x},\tilde{y})]
\int_{\mbb{N}^{-1}}[F(x,y+\xi,\tilde{x},\tilde{y}+\xi)-F(x,y,\tilde{x},\tilde{y})]
\,n(\d\xi)\cr
&~+[\gamma(x,y)-\gamma(\tilde{x},\tilde{y})]^+
\int_{\mbb{N}^{-1}}[F(x,y+\xi,\tilde{x},\tilde{y})-F(x,y,\tilde{x},\tilde{y})]
\,n(\d\xi)\cr
&~+[\gamma(x,y)-\gamma(\tilde{x},\tilde{y})]^-
\int_{\mbb{N}^{-1}}[F(x,y,\tilde{x},\tilde{y}+\xi)-F(x,y,\tilde{x},\tilde{y})]
\,n(\d\xi).
\end{split}
\end{equation}
Here and in what follows, $F'_x=\frac{\partial F(x,y,\tilde{x},\tilde{y})}{\partial x},
F''_{xx}=\frac{\partial^2 F(x,y,\tilde{x},\tilde{y})}{\partial x^2}$
and so on. Similarly, we can define $\tilde{A}F(x,y,\tilde{x},\tilde{y})$ for the case that $x<\tilde{x}$. By \eqref{defin coupling gener}, it is not hard to see that $\tilde{A}$ is indeed a coupling generator of $A$ defined by \eqref{generator}.

\btheorem
There exists a coupling $(Z(t),\tilde{Z}(t))_{t\geq0}=(X(t),Y(t),\tilde{X}(t),\tilde{Y}(t))_{t\geq0}$ whose generator is $\tilde{A}$.
\etheorem
\proof We consider the following SDE:
\beqlb\label{coupling proc}\left\{
\begin{aligned}
	& X(t) =X(0)-b\int_0^t X(s)\,\d s+\int_0^t\sqrt{2cX(s)}\,\d B(s)+\int_0^t\int_0^{X(s-)}\int_0^\infty \xi\,\tilde{M}(\d s,\d u,\d \xi), \cr
	&  Y(t)=Y(0)+\int_0^t\int_0^{Y(s-)}\int_0^{h(X(s-),Y(s-))}\int_{\mbb{N}^{-1}} \xi\,N(\d s,\d u,\d r,\d \xi),\cr
	& \tilde{X}(t)=\tilde{X}(0)-b\int_0^t \tilde{X}(s)\,\d s+\int_0^t\sqrt{2c\tilde{X}(s)}\,\d B^*(s)+\int_0^t\int_0^{\tilde{X}(s-)}\int_0^\infty \xi\,\tilde{M}(\d s,\d u,\d \xi), \cr
	& \tilde{Y}(t)=\tilde{Y}(0)+\int_0^t\int_0^{\tilde{Y}(s-)}\int_0^{h(\tilde{X}(s-),\tilde{Y}(s-))}\int_{\mbb{N}^{-1}} \xi\,N(\d s,\d u,\d r,\d \xi),
\end{aligned}
\right.
\eeqlb
where
\beqnn
B^*(t)=\left\{
\begin{aligned}
	& -B(t),\quad\qquad\quad ~~~t\leq T, \cr
	&  -2B(T)+B(t),\quad t>T
\end{aligned}
\right.
\eeqnn
and $T=\inf\{t>0:X(t)=\tilde{X}(t)\}.$ Clearly, $(B^*(t))_{t\geq0}$ is still a standard Brownian motion. By the results in Section 2, we can determine the unique strong solution $(Z(t),\tilde{Z}(t))_{t\geq0} := (X(t),Y(t),\tilde{X}(t),\tilde{Y}(t))_{t\geq0}$ to (\ref{coupling proc}). On the other hand, we can apply the $\mrm{It\hat{o}}$'s formula to the SDE (\ref{coupling proc}) to see that the infinitesimal generator of the process $(Z(t),\tilde{Z}(t))_{t\geq0}$ is indeed the coupling generator defined by (\ref{coupling gener}).\qed

Now let us define a function $F$ on $D\times D$ such that
\beqlb
F(x,y,\tilde{x},\tilde{y})=F(|x-\tilde{x}|+|y-\tilde{y}|)=
|x-\tilde{x}|+\theta|y-\tilde{y}|
\eeqlb
with $\theta>0$, the exact value will be determined later. It is easy to see that
\beqnn
F(|x-\tilde{x}|+|y-\tilde{y}|)\asymp
|x-\tilde{x}|+|y-\tilde{y}|.
\eeqnn
For our main result in this subsection, we require the following assumptions.

\begin{condition}\label{moment of n}
	$b > 0$ and $R_1:=\int_{\mbb{N}^{-1}}\xi\,n(\d\xi)\in(-\infty,0)$.
\end{condition}

\begin{condition}\label{upper bound of gamma}
For any $y\in\mbb{N}_+$, $y\mapsto m(y)y$
is non-decreasing and bounded in the sense of $R_2:=\sup\limits_{y\in\mbb{N}_+}m(y)y<\infty$.
\end{condition}

\begin{remark}
	
	(1). $b>0$ means that $X$ is a subcritical CB-process; see, e.g., \cite{L11,L20}. There is not any restrictions on L\'evy noises in this paper. It is reasonable since we mainly focus on the exponential ergodicity in the $L^1$-Wasserstein distance. When considering the exponential ergodicity in the total variation distance for CB-processes, some additional assumptions on L\'evy noises are needed; see \cite{LM15,FJKR23}. In the setting of state-dependent branching cases, one has to make some restrictions on L\'evy noises even in the $L^1$-Wasserstein distance; see \cite{LW20}. Furthermore, we mention that either in other distances $W_d$ or state-dependent cases, more complicated functions are needed.

(2). $R_1 < 0$ of Condition \ref{moment of n} actually means that the associated first moment of offerspring of each individual strictly less than 1, i.e. $\sum_{j}jp_j<1$. In this case, the process $Y$ is called the subcritical case of continuous-time Markov branching processes when $h$ is a positive constant; see, e.g., \cite[pp.112]{AN72}.
	
(3). Condition \ref{upper bound of gamma} holds when $m(\cdot)\equiv0$. In this case, $(Y(t))_{t\geq0}$ is a standard continuous time branching process with branching rate $r>0$ and offspring $(p_\xi,\xi\in\mbb{N})$. Condition \ref{upper bound of gamma} also holds when $m(y)=\frac{1}{y+1}$. In this case, as the number of cells increasing, the rate of cell division is getting slower.
\end{remark}
	
We now present the main result.

\begin{theorem}
	Suppose that Conditions \ref{moment of n}--\ref{upper bound of gamma} are satisfied. Then there are constants $\lambda_0>0$ such that for any $(x,y),(\tilde{x},\tilde{y})\in D$, (\ref{exponential contract}) holds.
\end{theorem}
\proof We shall first give some estimates of $\tilde{A}F
(|x-\tilde{x}|+|y-\tilde{y}|)$. By (\ref{coupling gener}) and Taylor's formula, we have
\beqnn
\ar\ar\tilde{A}F(|x-\tilde{x}|+|y-\tilde{y}|)\cr
\ar\ar\quad \leq-bx\frac{x-\tilde{x}}{|x-\tilde{x}|}
+b\tilde{x}\frac{x-\tilde{x}}{|x-\tilde{x}|}\cr
\ar\ar\qquad~+\theta[\gamma(x,y)-\gamma(\tilde{x},\tilde{y})]^+
\int_{\mbb{N}^{-1}}\Big[|y-\tilde{y}+\xi|-|y-\tilde{y}|\Big]
\,n(\d\xi)\cr
\ar\ar\qquad~+\theta[\gamma(x,y)-\gamma(\tilde{x},\tilde{y})]^-
\int_{\mbb{N}^{-1}}\Big[|y-\tilde{y}-\xi|-|y-\tilde{y}|\Big]
\,n(\d\xi)\cr
\ar\ar\quad\leq-b|x-\tilde{x}|+\theta[\gamma(x,y)-\gamma(\tilde{x},\tilde{y})]^+\mbf{1}_{\{y>\tilde{y}\}}
\int_{\mbb{N}^{-1}}\xi\,n(\d\xi)\cr
\ar\ar\qquad~+\theta[\gamma(x,y)-\gamma(\tilde{x},\tilde{y})]^+\mbf{1}_{\{y\leq\tilde{y}\}}
\int_{\mbb{N}^{-1}}|\xi|\,n(\d\xi)+\theta[\gamma(x,y)-\gamma(\tilde{x},\tilde{y})]^-\mbf{1}_{\{y>\tilde{y}\}}
\int_{\mbb{N}^{-1}}|\xi|\,n(\d\xi)\cr
\ar\ar\qquad~+\theta[\gamma(x,y)-\gamma(\tilde{x},\tilde{y})]^-\mbf{1}_{\{y\leq\tilde{y}\}}
\int_{\mbb{N}^{-1}}\xi\,n(\d\xi).
\eeqnn
Since $\gamma(x,y)=ry+xm(y)y$, we arrive at
\begin{equation*}
\begin{split}
\tilde{A}F(|x-\tilde{x}|+|y-\tilde{y}|)&\leq-b|x-\tilde{x}|-\theta n(\{-1\})(\gamma(x,y)-\gamma(\tilde{x},\tilde{y}))\frac{y-\tilde{y}}{|y-\tilde{y}|}\\
&~+\theta\int_{\mbb{N}}\xi\,n(\d\xi)|\gamma(x,y)-\gamma(\tilde{x},\tilde{y})|\\
&=-b|x-\tilde{x}|-\theta n(\{-1\})(ry-r\tilde{y}+xm(y)y-\tilde{x}m(\tilde{y})y)\frac{y-\tilde{y}}{|y-\tilde{y}|}\\
&~+\theta\int_{\mbb{N}}\xi\,n(\d\xi)|ry-r\tilde{y}+xm(y)y-\tilde{x}m(\tilde{y})y|.
\end{split}
\end{equation*}
Notice that $\int_{\mbb{N}^{-1}}|\xi|\,n(\d\xi)\in(0,\infty)$ by $R_1 \in (-\infty, 0)$ of Condition \ref{moment of n}. For the case of $x>\tilde{x}$, it follows from Condition \ref{upper bound of gamma} that
\beqnn
\tilde{A}F(|x-\tilde{x}|+|y-\tilde{y}|)\ar\leq\ar
-b|x-\tilde{x}|-\theta n(\{-1\})\Big(r|y-\tilde{y}|+(x-\tilde{x})m(y)y\frac{y-\tilde{y}}{|y-\tilde{y}|}\cr
\ar\ar~~~~~~~~~~~~~~~~~~~~~~~~~~~~~~~~~~~~~~~~
+\tilde{x}(m(y)y-m(\tilde{y})\tilde{y})\frac{y-\tilde{y}}{|y-\tilde{y}|}\Big)\cr
\ar\ar~+\theta\int_{\mbb{N}}\xi\,n(\d\xi)\Big(r|y-\tilde{y}|+(x-\tilde{x})m(y)y+\tilde{x}|m(y)y-m(\tilde{y})\tilde{y}|\Big)\cr
\ar\leq\ar\! -\Big(\!b\!-\!\theta\int_{\mbb{N}^{-1}}\!|\xi|\,n(\d\xi)m(y)y\!\Big)|x-\tilde{x}|\!+\!\theta r R_1|y-\tilde{y}|\!+\!\tilde{x}\theta R_1|m(y)y\!-\!m(\tilde{y})\tilde{y}|\cr
\ar\leq\ar-\Big(b-\theta R_2\int_{\mbb{N}^{-1}}|\xi|\,n(\d\xi)\Big)|x-\tilde{x}|+\theta r R_1|y-\tilde{y}|\cr
\ar\leq\ar-\lambda_1F(|x-\tilde{x}|+|y-\tilde{y}|)
\eeqnn
for some $\lambda_1>0$ by setting $\theta:=\theta_1=\frac{b}{2R_2\int_{\mbb{N}^{-1}}|\xi|\,n(\d\xi)}$. When $x\leq\tilde{x}$, similarly we have
\begin{equation*}
\begin{split}
\tilde{A}F(|x-\tilde{x}|+|y-\tilde{y}|)&\leq
-(b+\theta R_1R_2)|x-\tilde{x}|+\theta R_1\Big(r|y-\tilde{y}|+x|m(\tilde{y})\tilde{y}-m(y)y|\Big)\\
&\leq-\lambda_2F(|x-\tilde{x}|+|y-\tilde{y}|).
\end{split}
\end{equation*}
for some $\lambda_2>0$ by setting $\theta:=\theta_2=\frac{-b}{2R_1R_2}>0.$
In conclusion, let $\theta=\theta_1\wedge\theta_2$ and $\lambda=\lambda_1\wedge\lambda_2$. Following similar arguments in step 2 of the proof for \cite[Theorem 3.1]{LW19}, we obtain the desired result.
\qed

%\begin{acks}[Acknowledgments]
%	The authors would like to thank Prof. Wang Jian for his constructive comments that improved the 	quality of this paper.
%\end{acks}

%%%%%%%%%%%%%%%%%%%%%%%%%%%%%%%%%%%%%%%%%%%%%%
%% Funding information, if any,             %%
%% should be provided in the                %%
%% funding section.                         %%
%%%%%%%%%%%%%%%%%%%%%%%%%%%%%%%%%%%%%%%%%%%%%%
%\begin{funding}
	
%	The first author was supported in part by National Key R\&D Program of China (No. 2022YFA1006000) and China Postdoctoral Science Foundation (No. 2022M720735).
	
%	The second author was supported in part by Guangdong Basic and Applied Basic Research Foundation (No. 2022A1515110986), Guangdong Young Innovative Talents Project (No. 2022KQNCX105) and NSFC grant (No. 12271029).
	
%	The third author was supported in part by National Key R\&D Program of China grant (No. 2022YFA1006102) and NSFC grant (No. 11831010).
%\end{funding}


\begin{thebibliography}{00}
	\small
	
	\bibitem{AN72} Athreya, K.B. and Ney, P.E. (1972): {\it Branching Processes.} New York, Heidelberg, Berlin.
	
	\bibitem{BT11} Bansaye, V. and Tran V.C. (2011): Branching Feller diffusion for cell division with parasite infection. {\it ALEA, Lat. Am. J. Probab. Math. Stat.} {\bf 8}: 95--127.
	
	\bibitem{C05} Chen, M. (2005): \textit{ Eigenvalues, Inequalities, and Ergodic Theory},  Springer, New York.
	
	\bibitem{CL21} Chen, S. and Li, Z. (2021): Continuous time mixed state branching processes and stochastic equations. \textit{ Acta Math. Sci. (Engl. Ed.)} \textbf{41(5)}: 1445--1473.
	
\bibitem{CG95} Cranston, M. and  Greven, A. (1995): Coupling and harmonic functions in the case of continuous time Markov processes. {\it Stoch. Process. Appl.} {\bf60}: 261--286.

	\bibitem{DL06} Dawson, D.A. and Li, Z. (2006): Skew convolution semigroups and affine Markov processes. \textit{Ann. Probab.} \textbf{34}: 1103--1142.
	
	\bibitem{DL12} Dawson, D.A. and Li, Z. (2012): Stochastic equations, flows and measure-valued processes. \textit{Ann. Probab.} {\bf 40}: 813--857.
	
 \bibitem{EK86} Ethier, S.N. and Kurtz, T.G. (1986): {\it Markov Processes: Characterization and Convergence}. John Wiley and Sons, New York.

\bibitem{FJKR23} Friesen, M., Jin, P., Kremer, J. and R\"udiger, B. (2023): Exponential ergodicity for stochastic equations of nonnegative processes with jumps. {\it ALEA, Lat. Am. J. Probab. Math. Stat.} {\bf 20}: 593--627.

	\bibitem{FL10} Fu, Z. and Li, Z. (2010): Stochastic equations of non-negative processes with jumps. \textit{Stoch. Process. Appl.} \textbf{120}: 306--330.
	
	\bibitem{JX22} Ji, L. and Xiong, J. (2022): Superprocesses for the population of rabbits on grassland. \textit{ Proc. Steklov Inst. Math.} {\bf 316}(1): 195--208.
	
	\bibitem{KW71} Kawazu, K. and Watanabe, S. (1971): Branching processes with immigration and related limit theorems. \textit{Theory Probab. Appl.} {\bf 16}: 36--54.
	
\bibitem{L92} Lindvall, T. (1992): {\it Lectures on the coupling method.} Wiley Series in Probability and Mathematical Statistics: Probability and Mathematical Statistics, John Wiley \& Sons, Inc., New York.

\bibitem{LLWZ23+} Li, P., Li, Z., Wang, J. and Zhou, X. (2023+): Exponential ergodicity of branching processes with immigration and competition. To appear in {\it Annales de l'Institute Henri Poincar\'{e}-Probabilit\'{e}s et Statistiques}, available at arxiv: 2205.15499v1.

\bibitem{LW20} Li, P. and Wang, J. (2020): Exponential ergodicity for general continuous-state nonlinear branching processes. {\it Electron. J. Probab.} {\bf 25}: article no. 125, 1--25.

	\bibitem{LYZ19} Li, P., Yang, X. and Zhou, X. (2019): A general continuous-state nonlinear branching process. \textit{Ann. Appl. Probab.} \textbf{29}(4): 2523--2555.
	
	\bibitem{L06} Li, Z. (2006): A limit theorem for discrete Galton-Watson branching processes with immigration. \textit{J. Appl. Probab.} {\bf 43}: 289--295.
	
	\bibitem{L11} Li, Z. (2011): \textit{ Measure-Valued Branching Markov Processes.} Probab. and Its Appl. Springer.
	
	\bibitem{LM15} Li, Z. and Ma, C. (2015): Asymptotic properties of estimators in a Cox-Ingersoll-Ross model. {\it Stoch. Process. Appl.} {\bf 125}: 3196--3233.
	
	\bibitem{L20} Li, Z. (2020): \textit{ Continuous-State Branching Processes with Immigration.} In: Jiao, Y. (eds) From Probability to Finance. Mathematical Lectures from Peking University. Springer, Singapore.
	
	\bibitem{Liu22} Liu, J. (2022+): Scaling Limits of Controlled Branching Processes. Available at arXiv: 2204.06796.
	
	\bibitem{LW19} Luo, D. and Wang, J. (2019): Refined basic couplings and Wasserstein-type distances for SDEs with L\'evy noises. \textit{Stoc. Process. Appl.} \textbf{129}: 3129--3173.
	
	\bibitem{M09} Ma, C. (2009): A limit theorem of two-type Galton-Watson branching processes with immigration. \textit{Stat. Prob. Lett.} {\bf 79}: 1710--1716.
	
	\bibitem{MS20} Marguet, A. and Smadi, C. (2020+): Parasite infection in a cell population with deaths. Available at arXiv: 2010. 16070v1.
	
	
	\bibitem{OW20} Osorio, L. and Winter, A. (2020+): Two level branching model for virus population under cell division. Available at arXiv: 2004.14352.
	
	\bibitem{V09} Villani, C. (2009): \textit{ Optimal Transport, Old and New}. Springer, Berlin.
	
\end{thebibliography}
\end{document}